\documentclass[journal]{IEEEtran}
\pdfoutput=1
\usepackage{amsmath}
\usepackage{amssymb}
\usepackage{graphicx}
\usepackage{epsfig}
\usepackage{algorithm}
\usepackage{algorithmic}
\usepackage{caption}

\usepackage{float}
\usepackage{graphicx}
\usepackage{mathtools}
\usepackage{subfigure}

\begin{document}

\newtheorem{Thm}{Theorem}[section]
\newtheorem{Cor}[Thm]{Corollary}
\newtheorem{Lem}[Thm]{Lemma}
\newtheorem{Prop}[Thm]{Proposition}
\newtheorem{Def}[Thm]{Definition}
\newtheorem{rem}[Thm]{Remark}

\newtheorem{Assump}[Thm]{Assumption}

\DeclarePairedDelimiter{\norm}{ \lVert }{ \rVert }
\DeclarePairedDelimiter{\abs}{\lvert}{\rvert}
\DeclarePairedDelimiter{\inner}{\langle}{\rangle}


\title{\textbf{Newton-Step-Based Hard Thresholding Algorithms for Sparse Signal Recovery }}

\author{Nan Meng\thanks{Nan Meng is with the School of Mathematics, University of Birmingham,
Edgbaston, Birmingham B15 2TT,  UK (e-mail:
nxm563@bham.ac.uk). }   and  Yun-Bin Zhao, \emph{Member, IEEE}\thanks{Yun-Bin Zhao is with the School of Mathematics, University of Birmingham,
Edgbaston, Birmingham B15 2TT,  UK, and also with the Shenzhen Research Institute of Big Data  (e-mail:
y.zhao.2@bham.ac.uk). }
}

\date{}

\maketitle

\begin{abstract} Sparse signal recovery or compressed sensing   can be formulated as certain sparse optimization problems.  The classic optimization theory indicates  that the Newton-like method often has a numerical advantage over the gradient method for   nonlinear optimization problems.
In this paper, we propose the so-called Newton-step-based iterative hard thresholding (NSIHT) and the Newton-step-based hard thresholding pursuit (NSHTP) algorithms for sparse signal recovery and signal approximation. Different from the traditional iterative hard thresholding (IHT)  and hard thresholding pursuit (HTP),    the proposed algorithms  adopts the Newton-like search direction instead of the steepest descent direction.
 A theoretical analysis for the proposed algorithms is carried out, and  some sufficient conditions for the guaranteed success of sparse signal recovery via these algorithms are established. Our results are shown under the restricted isometry property which is one of the standard assumptions widely used in the field  of compressed sensing and signal approximation. The empirical results obtained from synthetic data recovery indicate that the proposed algorithms are efficient signal recovery methods.  The numerical stability of our algorithms in terms of the residual reduction is also investigated through simulations.

\end{abstract}


\section{Introduction}

 The compressed sensing (CS) theory was first introduced by Donoho, Cand\`es, Romberg, and Tao (see, e.g., \cite{candes2006compressive, candes2004robust, donoho2006compressed}).  It goes beyond the restriction of the classic Nyquist-Shannon sampling theory in signal processing \cite{claude1934communication}.
The compressed sensing together with its sparse optimization models has  widely been studied and applied in diverse areas such as the image processing  \cite{elad2010, foucart2013mathematical,  starck2010sparse}, data separation  \cite{eldar2012compressed}, pattern recognition \cite{patel2011sparse}, and wireless Network communication \cite{ han2013compressive}, to name a few.
The mathematical foundation for CS and sparse optimization problems can be found in such references as \cite{elad2010, eldar2012compressed, foucart2013mathematical, zhao2018sparse}.

The sparse optimization model is a backbone for the development of the theory and algorithms for CS and many other aspects of signal processing. One of the important optimization models for CS  is the following minimization with a sparsity constraint:
\begin{equation} \label{PPSS}
    \min_{x} \{ \|y-Ax\|_2^2 : ~  \|x\|_0 \le k \},
\end{equation}
where the function $\|x\|_0$ denotes the number of nonzero components of the vector $x,$  $A \in \mathbb{R}^{m \times n}$ ($m \ll n$) is a measurement matrix, and $y$ is the acquired information (called measurements) of the unknown signal to recover.
Depending on application environments, the CS or signal recovery model can take other  forms such as
\begin{equation} \label{PPSS-01} \min_{x} \{ \|x\|_0 : ~ y=Ax \} \end{equation}
 and the form
\begin{equation} \label{PPSS-02} \min_x \{ \|x\|_0:  ~ \|y-Ax\|_2^2 \le \eta \}.
  \end{equation}
The model (\ref{PPSS-01})  is  widely studied in noiseless settings, while the model (\ref{PPSS-02}) is more plausible from a practical viewpoint due to the existence of noises.
These sparse optimization problems are known to be   NP-hard in general \cite{N95}.
Clearly, the main difficulty for solving such a sparse optimization problem lies in the combinatorial nature of locating the support of the unknown signal.

Several approaches can be used to possibly solve the sparse optimization problems arising from signal recovery or  approximation.
For instance,   $\ell_1$-minimization and reweighted $\ell_1$-minimization  \cite{candes2008enhancing, chen2001atomic, zhao2012reweighted} are widely used to solve the problems (\ref{PPSS-01}) and (\ref{PPSS-02}).  The dual-density-based method  \cite{zhao2018sparse,zhao2016constructing} and other   heuristic algorithms such as the orthogonal matching pursuit (OMP) \cite{donoho2006sparse, tropp2007signal}, compressive sampling matching pursuit  (CoSaMP) \cite{needell2009cosamp}, and subspace pursuit  (SP) \cite{dai2009subspace}    are also useful methods for solving the problems (\ref{PPSS-01}) and (\ref{PPSS-02}).
 The hard-thresholding-type methods provide a simple way to ensure that the iterates generated by an iterative algorithm are feasible to the model (\ref{PPSS}).   Two commonly used thresholding techniques are the hard thresholding (e.g., \cite{blumensath2009iterative}-\cite{BFH16}, \cite{foucart2013mathematical},  \cite{FS19}) and soft thresholding (e.g., \cite{donoho1995noising, donoho1994ideal, elad2006simple, fornasier2008iterative,  herrity2006sparse}).
The iterative hard thresholding (IHT) combined with a pursuit step forms the so called hard thresholding pursuit (HTP) in \cite{foucart2011hard}. Some latest developments and applications of IHT and HTP can be found in \cite{BKS19, DSZ16, JLZ18, KK18, LB19, SL17, ZXLD19}.
It is worth mentioning that Zhao \cite{zhao2019optimal} (see also in \cite{Zhaoluo2020})  proposed a new thresholding technique called the optimal $k$-thresholding. This technique  can successfully overcome the numerical oscillation phenomenon in existing HTP-type methods, and was shown to be an efficient thresholding technique from both theoretical and practical viewpoints.

The iterative scheme of the IHT algorithm
\[ x^{p+1} = {\cal H}_k(x^p + \lambda A^T(y-Ax^p))  \] was derived by adopting the steepest descent direction of the objective function of the problem (\ref{PPSS}) at the  iterate $x^p $ (see, e.g.,  \cite{blumensath2009iterative}-\cite{blumensath2010normalized}). The   ${\cal H}_k(\cdot)$ is called a hard thresholding operator which retains the largest $k$ magnitudes and zeroes out the rest components of a vector. The theoretical properties of IHT and its variants have been widely studied over the past decade (see, e.g.,  \cite{blumensath2009iterative}-\cite{blumensath2010normalized},   \cite{foucart2011hard}-\cite{foucart2013mathematical}).
The classic optimization theory  \cite{BV04, nocedal2006numerical} has shown that the Newton-like method usually has a faster convergence  than the gradient method, and thus it is generally more efficient than the steepest descent method for solving nonlinear optimization problems.
This motivates us to develop a new hard thresholding method which adopts a Newton-like search direction instead of the steepest descent direction.

Given a twice continuously differentiable function $f$ with gradient $\nabla f(x),$  when its Hessian matrix $\nabla^2 f(x) $ is nonsingular at the iterate $x^p,$ the classic Newton's method takes the following iterative scheme:
$$  x^{p+1} = x^p - \lambda (\nabla^2 f(x^p))^{-1} \nabla f(x^p),   $$
where $\lambda>0 $ is a certain stepsize that can be determined by a certain line search method (see \cite{nocedal2006numerical} for details).
 In compressed sensing scenarios, however, the Hessian matrix of the objective function in (\ref{PPSS}) is singular.  Although a direct use of the aforementioned Newton's method is difficult, one is still able  to develop a   Newton-like method  for the signal recovery problem (\ref{PPSS}). We describe the algorithms in detail in the next section.  The main work of this paper includes:  (a) We establish some sufficient conditions for the guaranteed success (i.e., convergence) of the proposed  algorithms for sparse signal recovery; (b) We investigate the numerical behaviour (stability and signal recovery capability) of the proposed algorithms through simulations.

The paper is organized as follows.  The algorithms are introduced in Section II. The convergence of the algorithms in noiseless situations is shown in Section III. The analysis of the algorithms in noisy settings is carried out in Section IV. The numerical results are discussed in Section V.

\section{Newton-Step-Based Hard Thresholding Algorithms}

We first introduce some notations.  We use  $A \in \mathbb{R}^{m \times n}$ with $m \ll n$ to denote the measurement matrix and  $y \in \mathbb{R}^m$  the measurements of the unknown signal $x^*\in \mathbb{R}^n. $
The $\ell_2$-norm is defined as $\|x\|_2 := \sqrt{\sum_i^n x_i^2}$, where $x \in \mathbb{R}^n$. Throughout the paper,  we use  $[N]$ to denote the set $\{ 1,\dots,n \}. $
Given a set $S \subset [N], $ $\overline{S} := [N] \setminus S $ is the complement set of $S. $
We use $ |S|=\textrm{card}(S)  $ to denote  the cardinality of the set $S$.
The support of a vector $x$ is denoted by $\textrm{supp} (x) := \{i: x_i \neq 0 \}$.
We use $I$ to denote the identity matrix.
Given a matrix $A, $  $A^T$ denotes its transpose, and the $\ell_2$-norm of  $A$ is defined as $\| A \|_{2} \coloneqq \sqrt{\lambda_{\max} (A^TA)}$.
For a given set $\Omega \subseteq [N]$, unless otherwise stated $x_\Omega$ is the subvector of $x$ with entries indexed by $\Omega$, and  $A_\Omega$ denotes a submatrix of the matrix $A$ with columns indexed by $\Omega. $

For convenience of discussion, we may write the problem (\ref{PPSS}) equivalently as
$$
    \min_{x} \{ \frac{1}{2} \|y-Ax\|_2^2 : ~ \|x\|_0 \le k \}.
$$
Let $f(x) = \frac{1}{2}\|y-Ax\|_2^2$.
The gradient and Hessian of $f(x)$ are given, respectively, by
\[ \nabla f(x) = -A^T(y-Ax),  ~~ \nabla^2 f(x) = A^TA . \]
Clearly, the matrix $A^TA$ is singular  since $A \in \mathbb{R}^{m \times n}$ and  $m \ll n$.
In order to develop a Newton-like method for the  model (\ref{PPSS}), we need to introduce a suitable modification of the Hessian of $f(x).$ A simple idea is to perturb the Hessian with a parameter $\epsilon > 0$  such that
$   A^TA + \epsilon I$ is positive definite,
where $I \in \mathbb{R}^{n \times n}$ is the identity matrix.
Such a perturbation of the Hessian leads to the following Newton-like iterative method for minimizing the unconstrained function $f(x) : $
\begin{equation} \label{NEWTON}  x^{new} = x^p +  \lambda (A^TA + \epsilon I)^{-1} A^T (y-Ax^p) ), \end{equation}
 where $ x^p $ is the current iterate, and $x^{new} $ is the new iterate generated as above.  However, the vector $ x^{new}$ may not be $k$-sparse and thus may not satisfy the constraint of the problem (\ref{PPSS}).
Therefore, to develop an iterative method that generates iterates satisfying the constraint of (\ref{PPSS}), it makes sense to consider the following iterative scheme:
\[ x^{p+1} = {\cal H}_k \left(x^p +  \lambda (A^TA + \epsilon I)^{-1} A^T (y-Ax^p) \right), \]
where the hard thresholding operator retains only the largest $k$ magnitudes of the vector $x^ {new} $ obtained by the Newton-like method (\ref{NEWTON}).
We now formatively state the algorithms for the problem (\ref{PPSS})  as follows.

\begin{algorithm} \caption{~ [Newton-Step-Based Iterative Hard Thresholding (NSIHT)]}
\begin{itemize} \label{IHTwithnewtonstep}
    \item  Input:   measurement matrix $A$, measurement vector $y$, sparsity level $k$, parameter $\epsilon > 0,$ and stepsize $\lambda$.
    \item  Iteration:
            \[ x^{p+1} = {\cal H}_k \left(x^p + \lambda (A^TA + \epsilon I)^{-1} A^T (y-Ax^p) \right). \]
    \item  Output:  The $k$-sparse vector $x^*$.
\end{itemize}
\end{algorithm}

It is possible to further reduce the  objective of the problem (\ref{PPSS}) by performing a pursuit step  as  traditional HTP algorithms, leading to the following algorithm called NSHTP.

\begin{algorithm}
\caption{ [Newton-Step-Based Hard Thresholding Pursuit (NSHTP)]}
\begin{itemize} \label{HTPwithnewtonstep}
    \item  Input:  measurement matrix $A$, measurement vector $y$, sparsity level $k$, parameter $\epsilon > 0,$ and stepsize $\lambda$.
    \item  Iteration:
            $$ \overline{x}^{p}  = {\cal H}_k \left(x^p + \lambda (A^TA + \epsilon I)^{-1} A^T (y-Ax^p) \right), $$
          \begin{equation} \label{PS}  x^{p+1}  = \arg \min_z \{\|y-Az\|_2^2: ~ \textrm{supp} (z) \subseteq \textrm{supp} (\overline{x}^p) \}.
\end{equation}
    \item  Output:  The $k$-sparse vector $x^*$.
\end{itemize}
\end{algorithm}
The  step (\ref{PS}) is called a pursuit step, which is to  minimize the objective function over the support of the iterate generated by the NSIHT.
As a result, $\|y-Ax^{p+1} \|_2 \leq \| y-A  \overline{x}^{p} \|_2.$ The solution to the pursuit step (\ref{PS}) can have a closed form when the sparsity level $k$ is low.  In fact, any $k$ columns of $A$ are linearly independent when $k$ is low enough, for instance, lower than the spark of the matrix. Denote the support set of $ \overline{x}^p$ by $\Omega.$ Then the problem (\ref{PS}) becomes an unconstrained minimization problem   $ \min \|y - A_\Omega z_\Omega\|_2, $  to which when $A_\Omega$ has linearly independent columns,
the  solution to the   problem  is given explicitly  as $x_\Omega = ((A_\Omega)^T A_\Omega)^{-1} (A_\Omega)^T y .$

In later sections, we carry out theoretical analysis for the above algorithms and establish the sufficient conditions for the guaranteed success of   signal recovery via these algorithms in both noiseless and noisy environments. Before doing so, we first introduce the restricted isometry constant (RIC) of a given matrix, which is very useful tool and has been widely used in the  CS literature.

\begin{Def}
 \cite{candes2005decoding}
The $q$-th order restricted isometry constant (RIC) $\delta_q$ of a matrix $A \in \mathbb{R}^{m \times n}$ is the smallest number $\delta_q \in (0,1]$ such that
\[ (1-\delta_q)\|x\|_2^2\le \|Ax\|_2^2 \le (1+\delta_q)\|x\|_2^2 \]
for all $q$-sparse vectors $x$.
\end{Def}

In the above definition, $q$ is an integer number.  A vector $x$ is called $q$-sparse if the number of nonzero components of $x$ is smaller than or equal to  $q$, i.e., $\norm{x}_0 \le q$.
An alternative expression of the $q$-th order restricted isometry constant (see \cite{foucart2011hard, foucart2012sparse}) is $\delta_q   =\max_{S\subset [N],|S| \le q}\|A^T_SA_S-I\|_{2}. $ In this paper,
the  matrix $A$ is said to satisfy the $q$-th order restricted isometry property (RIP) if $\delta_q$ is  smaller than 1.  It was shown that the random matrices such as Bernoulli and Gaussian random matrices may satisfy the RIP with a dominant probability \cite{candes2004robust, candes2005decoding, foucart2013mathematical, rauhut2010compressive}.
The matrix satisfying the RIP is often called a RIP matrix, which is widely used in the analysis of  various compressed sensing algorithms.

\section{Guaranteed success of NSIHT in noiseless settings}
We first analyze the NSIHT algorithm when the measurements $y:=Ax$ of the signal $x$ are accurate. Without loss of generality, we assume in this section that the target signal  $x$ is sparse. Our analysis here can be easily extended to the case when the signal can be sparsely approximated.
We show that the success of   signal recovery can be guaranteed under a standard RIP assumption and the suitable choice of the parameter $\epsilon$ as well as the stepsize $\lambda.$
Before showing the main result  of this section, we need  some useful technical results.

\begin{Lem}\textup{\cite{foucart2013mathematical}} \label{lemma1}
Let $ A^{m\times n}$ with $ m\ll n$ be a measurement matrix.  Given vectors $u,v \in \mathbb{R}^{n}, $
if $|\text{supp}(u) \cup \text{supp}(v)| \le t, $ one has
\[ \abs{\inner{u,(I-A^TA)v}} \le \delta_t \|u\|_2 \|v\|_2 . \]
\end{Lem}

The next lemma is key for our later analysis.

\begin{Lem} \label{lemma2}   Let $ A^{m\times n}$ with $ m\ll n$ be a measurement matrix.
Given vectors $u,v \in \mathbb{R}^{n}$ and an index set $\Omega \subset{[N]},$ if   $ \sigma^2_1 < \epsilon$ and  $\lambda \le \epsilon + \sigma_m^2,$ where $\sigma_1,\sigma_m$ are  the largest and smallest singular values of the matrix $A,$ respectively, then one has
 \begin{align*}  & | \langle u,    (I - \lambda (A^TA + \epsilon I)^{-1} A^TA) v \rangle  |  \\
              & \le ( \delta_t + \sigma_1^2-\frac{\lambda\sigma_1^2}{\epsilon+\sigma_1^2} )
    \|u\|_2 \|v\|_2 ~~ \textrm{ if }   |\textrm{supp}(u) \cup \textrm{supp}(v)| \le t,
\end{align*}
and
 \begin{align*}   & \|[ (I -  \lambda (A^TA+\epsilon I)^{-1} A^TA ) v ]_\Omega \|_2  \\
 & \le ( \delta_t + \sigma_1^2-\frac{\lambda\sigma_1^2}{\epsilon+\sigma_1^2} ) \|v\|_2  ~~ \textrm{ if }  |\Omega \cup \textrm{supp}(v)| \le t .
\end{align*}

\end{Lem}

\emph{Proof.}
The matrix $(\epsilon I +A^TA)^{-1}$ can be written as
\begin{equation} \label{star}
    (\epsilon I +A^TA)^{-1} = \frac{1}{\epsilon} (I + \frac{1}{\epsilon} A^TA)^{-1}.
\end{equation}
It is well known that if a square matrix $M$ satisfies $\|M\|_2<1$, then
\begin{equation} \label{matrixexpansion}
    (I+M)^{-1} = I-M + M^2 - M^3 + \dots = I + \sum_{i=1}^\infty (-1)^i M^i.
\end{equation}
In order to use (\ref{matrixexpansion}), we choose the parameter $ \epsilon $ such that
\[ \|\frac{1}{\epsilon} A^TA \|_2 = \sqrt{\frac{1}{\epsilon} \lambda_{\max}( A^TA)} <1, \]
i.e.,  $ \epsilon > \sigma_1^2,$ where $\lambda_{\max} (A^TA) = \sigma^2_1 $ is the largest eigenvalue of $ A^TA, $
Under such a choice of $ \epsilon,$  by (\ref{matrixexpansion}), the matrix $(I + \frac{1}{\epsilon} A^TA)^{-1}$ can be expanded as
\begin{align}  \label{starstar}
    (I+\frac{1}{\epsilon} A^TA)^{-1} & =I + \sum_{i=1}^\infty (-1)^i \left(\frac{1}{\epsilon} A^TA \right) ^i  \nonumber \\
     & = I + \sum_{i=1}^\infty (-\frac{1}{\epsilon})^i (A^TA) ^i.
\end{align}
By using (\ref{star}) and (\ref{starstar}), we have
\begin{equation*} \begin{aligned}
      I - \lambda   &   (A^TA+\epsilon I)^{-1} A^TA \\
       & = I - \frac{\lambda}{\epsilon} (I + \frac{1}{\epsilon} A^TA)^{-1} A^TA \\
      & = I - \frac{\lambda}{\epsilon} \left[  I + \sum_{i=1}^\infty (-\frac{1}{\epsilon})^i (A^TA) ^i  \right] A^TA \\
      & =   I +   \lambda \left[ -\frac{1}{\epsilon} A^TA  + \sum_{i=1}^\infty (-\frac{1}{\epsilon})^{i+1} (A^TA) ^{i+1}  \right]  \\
      &   =   I +   \lambda    \sum_{j=1}^\infty (-\frac{1}{\epsilon})^j  (A^TA) ^j  \\
    & = ( I - A^TA )  +  B  ,
\end{aligned} \end{equation*}
where
$$
    B  = A^TA + \lambda    \sum_{j=1}^\infty (-\frac{1}{\epsilon})^j  (A^TA) ^j .
$$
Then
  \begin{align} \label{addinner}
   | \langle u,  &  (I - \lambda (A^TA + \epsilon I)^{-1} A^TA) v \rangle  | \nonumber \\
&  = | \langle u, ( I - A^TA + B) v \rangle |  \nonumber \\
&  = | \inner*{u,(I - A^TA)v} + \inner*{u, B v} |  \nonumber \\
&  \le | \inner*{u,(I - A^TA)v} | + | \inner*{u, B v} | ,
\end{align}
where the last relation follows from the triangle inequality.
By Lemma \ref{lemma1}, the first term of the right-hand side of (\ref{addinner}) can be bounded as
\begin{equation} \label{firstone}
    \left| \inner*{u,(I - A^TA)v} \right| \le \delta_{t} \norm{u}_2 \norm{v}_2 ,
\end{equation} where $ t= |\textrm{supp} (u) \cup \textrm{supp} (v)|. $
We now bound the   term  $|\inner*{u, B v} |. $ Note that
\begin{equation} \label{inner2}
    \left| \inner*{u, B v} \right| \le \norm*{B}_2 \norm{u}_2 \norm{v}_2.
\end{equation}
It is sufficient to bound the  $\norm*{B}_2.$  To this need,  let  $ A=U \Sigma V^T $
 be the singular value decomposition, where $U,V$ are two orthogonal matrices with $U\in \mathbb{R}^{m \times m}$ and $V \in \mathbb{R}^{n\times n}$ and
$$ \Sigma = \left[ \begin{array}{cccc|c} \sigma_1 &&& &\\
                    & \sigma_2 && & 0_{m\times (n-m)}\\
                    & & \ddots & & \\
                    &&& \sigma_m &\end{array} \right]_{m\times n} ,$$
where $\sigma_i,$ $i = 1, \dots , m$ are singular values of $A$ satisfying $\sigma_1 \ge \sigma_2 \ge \dots \ge \sigma_m$. Under the choice of $ \epsilon,$  it is very easy to see that
the eigenvalues of the matrix $B$ are given as
\[ \begin{aligned}
    & \text{eig}(B)_i = \sigma_i^2-\frac{\lambda\sigma_i^2}{\epsilon+\sigma_i^2}, && i=1,\dots,m, \\
    & \text{eig}(B)_i = 0, && i=m+1,\dots, n.
\end{aligned} \]
Also, it is easy to verify that if $\lambda$ is chosen such that $\lambda \le \epsilon + \sigma_m^2$, where $\sigma_m$ is the smallest singular value of the matrix $A$, then  $\sigma_1^2-\frac{\lambda\sigma_1^2}{\epsilon+\sigma_1^2}$ is the largest eigenvalue of $B, $ i.e.,
\begin{equation} \label{norm}
    \norm*{B}_2 = \sigma_1^2-\frac{\lambda\sigma_1^2}{\epsilon+\sigma_1^2} .
\end{equation}
Substituting (\ref{norm}) into (\ref{inner2}) leads to
\begin{equation} \label{secondone}
    \left| \langle u, B v \rangle \right| \le  \left( \sigma_1^2-\frac{\lambda\sigma_1^2}{\epsilon+\sigma_1^2} \right) \norm{u}_2 \norm{v}_2 .
\end{equation}
Combining  (\ref{addinner}), (\ref{firstone}) and (\ref{secondone})  immediately yields the first inequality of the lemma. The second inequality in the lemma follows immediately from the first one.  In fact,  for any set $\Omega \subset [N]$,  we see that
\[ \begin{aligned}
    & \left\| \left[ (I - \lambda (A^TA+\epsilon I)^{-1} A^TA ) v \right]_\Omega\right\|_2^2 \\
    & =
   \langle  [ (I - \lambda (A^TA+\epsilon I)^{-1} A^TA ) v  ]_\Omega ,  \\
     & ~~~~~~~~~~~~~  [ (I - \lambda (A^TA+\epsilon I)^{-1} A^TA ) v  ]_\Omega \rangle \\
    & = \langle  [ (I - \lambda (A^TA+\epsilon I)^{-1} A^TA ) v  ]_\Omega ,  \\
     & ~~~~~~~~~~~~~  (I - \lambda (A^TA+\epsilon I)^{-1} A^TA ) v \rangle \\
    & \le ( \delta_t + \sigma_1^2 - \frac{\lambda \sigma_1^2}{\epsilon + \sigma_1^2} ) \\
     & ~~~~~~~~~~  \cdot \| [ (I - \lambda (A^TA+\epsilon I)^{-1} A^TA ) v  ]_{\Omega}  \|_2 \|v\|_2 ,
\end{aligned} \]
where the final inequality follows from the first inequality in the lemma that we just shown above by setting $u= (I - \lambda (A^TA+\epsilon I)^{-1} A^TA ) v  ]_\Omega. $  Dividing through the above inequality by $\| \left[ (I - \lambda (A^TA+\epsilon I)^{-1} A^TA ) v \right]_\Omega\|_2$ leads to the second inequality of the lemma.  ~~~ $ \Box$

The following property for the hard thresholding operator ${\cal H}_k$ is implied from the analysis of IHT by Foucart \cite{foucart2011hard}. However, we include a simple proof here for completeness.

\begin{Lem} \label{lem33}
For any vector $z \in \mathbb{R}^n $ and  any $k$-sparse vector $x\in \mathbb{R}^n$, one has
\[ \| {\cal H}_k(z) - x \|_2 \le \sqrt{3} \| (z - x)_{S \cup S^{*}} \|_2 \]
where  $S = \textrm{supp}(x)$ and  $S^{*} = \textrm{supp} ( {\cal H}_k(z) ). $
\end{Lem}

\emph{Proof.} Let $S$ and $S^*$ be defined as above.
Since $S^*$ is the support for the largest $k$ magnitudes  of $z,$  we have $\| z_S \|_2^2 \le \| z_{S^*} \|_2^2 . $ Eliminating those elements in $S \cap S^*$ leads to
\[  \| z_{S \setminus S^*} \|_2 \le \| z_{S^* \setminus S} \|_2 = \| (z - x)_{S^* \setminus S} \|_2, \]
where the equality follows from the fact $x_{S^* \setminus S}=0.$
Since $[{\cal H}_k(z)]_{S \setminus S^*}=0$,  we also have
\begin{align*}
    \| z_{S \setminus S^*} \|_2  & = \| (x - {\cal H}_k(z) - x + z)_{S \setminus S^*} \|_2 \nonumber \\
    &  \ge \| (x - {\cal H}_k(z))_{S \setminus S^*} \|_2 - \| (z-x)_{S \setminus S^*} \|_2 .
\end{align*}
Combing the above two relations, we get
\begin{align*}
    \| (x - {\cal H}_k(z))_{S \setminus S^*} \|_2
    &  \le  \| (z-x)_{S \setminus S^*} \|_2 + \| (z-x)_{S^* \setminus S} \|_2 \\
     &   \le  \sqrt{2} \| (x-z)_{S \cup S^*} \|_2 .
 \end{align*}
Therefore,
\[ \begin{aligned}
    \| x-{\cal H}_k(z) \|_2^2 & = \| \left(x-{\cal H}_k(z)\right)_{S^*} \|_2^2 + \| \left(x-{\cal H}_k(z)\right)_{\overline{S^*}} \|_2^2 \\
    & = \| \left(x-z\right)_{S^*} \|_2^2 + \| \left(x-{\cal H}_k(z)\right)_{S \setminus S^*} \|_2^2 \\
 &  \le \| \left(x-z\right)_{S^*} \|_2^2 + 2 \| (x-z)_{S \cup S^*} \|_2^2 \\
 &   \le 3 \| \left(x-z\right)_{S \cup S^*} \|_2^2,
\end{aligned} \]
as desired.  ~~~ $ \Box $

The main result of this section is summarized as follows.

\begin{Thm} \label{Thm1}
Let the restricted isometry constant $\delta_{3k}$ of the measurement matrix $A\in \mathbb{R}^{m \times n}$ satisfies that
$ \delta_{3k} < \frac{1}{\sqrt{3}} .  $
Let $\epsilon$ be a given parameter  satisfying
\begin{equation} \label{cond-epsilon} \epsilon > \max \left\{ \sigma_1^2, ~   \left( \frac{\sigma_1^2-\sigma_m^2}{\frac{1}{\sqrt{3}} -\delta_{3k}} - 1 \right) \sigma_1^2  \right\}. \end{equation}
If the stepsize $\lambda$ in NSIHT is chosen such that
\begin{equation} \label{rangeoflambda1}
  \epsilon + \sigma_1^2 - \left( \frac{1}{\sqrt{3}}-\delta_{3k}\right) \frac{\epsilon+\sigma_1^2}{\sigma_1^2}  <\lambda \leq \epsilon +\sigma_m^2  ,
\end{equation}
where $\sigma_1$  and $\sigma_m$ are the largest and smallest singular values of $A$, respectively,
then for any $k$-sparse signal $x\in\mathbb{R}^n$  with accurate measurements $y: =Ax, $  the sequence $\{x^n\}$ generated by the NSIHT algorithm converges to the signal $x$.
\end{Thm}

\emph{Proof.}
It is sufficient to prove that there exists a constant $\rho < 1$ such that the sequence generated by the NSIHT satisfies the relation below:
    \[ \|x^{p+1}-x\|_2 \le \rho \|x^p -x\|_2 , \]
which ensures the convergence of the  $\{x^p \} $ to $x.$  Note that  $ y=A x.$ We define
\begin{align*}  u^p    : = & x^p + \lambda (A^TA + \epsilon I)^{-1} A^T (y-Ax^p)  \nonumber  \\
  =  & x^p + \lambda (A^TA + \epsilon I)^{-1} A^TA (x-x^p) . \label{unexpression}
\end{align*}
Denote by $S = \text{supp}(x)$ and $S^{p+1} = \text{supp}(x^{p+1})$.
By the structure of the algorithm,  $ x^{p+1} = {\cal H}_k(u^{p}). $  We immediately have the following relation from  Lemma \ref{lem33}:
\begin{equation} \label{shizi3}
    \|x^{p+1} - x\|_2
    \le \sqrt{3} \| ( u^p - x )_{ S^{p+1} \cup S }\|_2 .
\end{equation}
From
\[ u^p - x = x^p - x + \lambda (A^T A + \epsilon I)^{-1} A^T A (x - x^p) ,  \]
we immediately have that
\begin{align}  \label{shizi66}
   &  \| (u^p  - x)_{ S^{p+1} \cup S } \|_2 \nonumber \\
    & = \|[ \left( I - \lambda (A^TA + \epsilon I)^{-1} A^T A \right) (x^p - x) ]_{ S^{p+1} \cup S } \|_2 .
\end{align}
Since $|S^{p+1}| \le k$, $|S| \le k$ and $|\textrm{supp}(x^p)| \le k$, we easily see that
\[  |S^{p+1} \cup S \cup \textrm{supp} (x^p - x)| \le 3k . \]
Under the condition $ \epsilon > \sigma_1^2,$ it follows from  Lemma \ref{lemma2}  that
\begin{align} \label{shizi55}
  &   \left\|\left[ \left( I - \lambda (A^TA + \epsilon I)^{-1} A^T A \right) (x^p - x) \right]_{ S^{p+1} \cup S } \right\|_2  \nonumber \\
    & \le
   \left ( \delta_{3k} + \sigma_1^2 - \frac{\lambda\sigma_1^2}{\epsilon+\sigma_1^2} \right)
\|x^p - x\|_2 ,
\end{align}
where $\lambda \le \epsilon + \sigma_m^2. $
Merging (\ref{shizi3}), (\ref{shizi66}) and (\ref{shizi55}),  we obtain
\begin{equation}  \label{condition1} \|x^{p+1} - x\|_2 \le \sqrt{3}
\left( \delta_{3k} + \sigma_1^2-\frac{\lambda\sigma_1^2}{\epsilon+\sigma_1^2} \right )
\|x^p - x\|_2  ,
\end{equation}
where the constant
$$
    \delta_{3k} + \sigma_1^2-\frac{\lambda\sigma_1^2}{\epsilon+\sigma_1^2} > 0
$$
  is guaranteed under the condition $\lambda \leq  \epsilon +\sigma_m^2$.
In addition, to ensure the convergence of the iterates $\{x^p\},$ the constant coefficient of the right-hand side of (\ref{condition1}) must be smaller than 1, i.e.,
$$
    \rho := \sqrt{3}
\left( \delta_{3k} + \sigma_1^2-\frac{\lambda\sigma_1^2}{\epsilon+\sigma_1^2} \right ) < 1 ,
$$
which is guaranteed if $\lambda $ is taken such that
$$
    \lambda > \epsilon + \sigma_1^2 - \left( \frac{1}{\sqrt{3}}-\delta_{3k}\right) \frac{\epsilon+\sigma_1^2}{\sigma_1^2}  .
$$
So the  range of $\lambda$ is determined as
\begin{equation} \label{lambda}
    \epsilon + \sigma_1^2 - \left( \frac{1}{\sqrt{3}}-\delta_{3k}\right) \frac{\epsilon+\sigma_1^2}{\sigma_1^2} <   \lambda  \leq \epsilon +\sigma_m^2.
\end{equation}
To make sure the existence of such a range, it suffices to require that
\[ \epsilon + \sigma_m^2 >  \epsilon + \sigma_1^2 - \left( \frac{1}{\sqrt{3}}-\delta_{3k}\right) \frac{\epsilon+\sigma_1^2}{\sigma_1^2}  . \]
This is guaranteed by the following choice of $\epsilon: $
\[ \epsilon > \left ( \frac{\sigma_1^2-\sigma_m^2}{\frac{1}{\sqrt{3}} -\delta_{3k}} - 1 \right) \sigma_1^2 , \]
which  together with $\epsilon>\sigma_1^2 $ leads to the condition (\ref{cond-epsilon}).     Such a choice of $\epsilon$ ensures that if the stepsize $\lambda$ is chosen as (\ref{lambda}), then the sequence generated by NSIHT converges to the target vector $x$ under the condition $\delta_{3k} < \frac{1}{\sqrt{3}}. $  ~~ $ \Box $\\

   For the special choice of the stepsize $\lambda=\epsilon$, following the same proof above, it can be seen that to ensure the convergence of the algorithm, $\epsilon$ should be chosen such that the following inequality is guaranteed:
\[ \sigma_1^2 - (\frac{1}{\sqrt{3}}-\delta_{3k}) \frac{\epsilon+\sigma_1^2}{\sigma_1^2} < 0 , \]
which  implies that
\[ \epsilon >\left ( \frac{\sigma_1^2}{\frac{1}{\sqrt{3}} -\delta_{3k}} - 1 \right) \sigma_1^2 . \]
Thus, together with the requirement $\epsilon > \sigma^2_1,$  the parameter $\epsilon$ is chosen to  satisfy the following condition in order to ensure the convergence of the  NSIHT with $\lambda=\epsilon:$
\begin{equation} \label{epsilon2}
    \epsilon > \max \left\{ \sigma_1^2 , \left( \frac{\sigma_1^2}{\frac{1}{\sqrt{3}} -\delta_{3k}} - 1 \right) \sigma_1^2 \right\} .
\end{equation} We summarize this result  as follows.

\begin{Cor}
Let the restricted isometry constant $\delta_{3k}$ of the matrix $A\in \mathbb{R}^{m \times n}$ satisfies that
$ \delta_{3k} < \frac{1}{\sqrt{3}},$ and let  $\epsilon$ be a given parameter  satisfying
(\ref{epsilon2}), where $\sigma_1$  is the largest singular value of $A$.
Then for any $k$-sparse signal $x\in\mathbb{R}^n$  with accurate measurements $y: =Ax, $  the sequence $\{x^n\}$ generated by the NSIHT algorithm with $\lambda =\epsilon $ converges to  $x. $
\end{Cor}

\section{Analysis of NSIHT and NSHTP in noisy settings}
In this section, we analyze both  NSIHT and NSHTP for signal recovery in noise scenarios, where  the signal may not necessarily  be $k$-sparse and the measurements $y=Ax+e$   is inaccurate with error $e.$  The recovery error bounds can be established for such situations under the  same or similar  conditions of Theorem \ref{Thm1}.
We first make an observation for the norm of the  matrix  $A (A^TA+\epsilon I)^{-1}$ that will be used in our analysis.
Applying singular value decomposition $A = U \Sigma V^T,$ we immediately see that
\[ \begin{aligned}  A (A^TA+\epsilon I)^{-1}
        = & U \begin{bmatrix}
        \frac{\sigma_1}{\epsilon + \sigma_1^2} \\
          &  \ddots &  &  0_{m \times (n-m)} \\
         &  & \frac{\sigma_m}{\epsilon + \sigma_m^2} \\
    \end{bmatrix} V^T,
\end{aligned} \]
where $U \in \mathbb{R}^{m \times m}$ and $V \in \mathbb{R}^{n \times n}$ are orthogonal matrices.
Clearly, we have
\[ \frac{\sigma_1}{\epsilon + \sigma_1^2} \ge \frac{\sigma_2}{\epsilon + \sigma_2^2} \ge \dots \ge \frac{\sigma_m}{\epsilon + \sigma_m^2} , \]
which implies that
\begin{equation} \label{xingxing2}
    \|A (A^TA+\epsilon I)^{-1}\| _2 = \frac{\sigma_1}{\epsilon + \sigma_1^2}.
\end{equation}

In the remainder of this section, we use $x_S\in \mathbb{R}^n$ to denote the vector obtained from $x\in \mathbb{R}^n $ by retaining the $k$ largest magnitudes of $x$ and zeroing out its other entries. Thus, in this section, the vector $x_S$ has the same dimension of $x$, and we denote by $x_{\overline{S}}= x-x_S.$

\begin{Thm}\label{thm2}
Let $\sigma_1$ and $\sigma_m$ be the largest and the smallest singular values of $A\in \mathbb{R}^{m \times n}. $
Suppose that the restricted isometry constant of $A$ satisfies that
$ \delta_{3k} < \frac{1}{\sqrt{3}}.  $
Let $y \coloneqq Ax+e$ be the measurements of the signal $x\in\mathbb{R}^n$ and $e$ is the measurement error. Let $\epsilon$ be  a given parameter satisfying (\ref{cond-epsilon}).
If the stepsize  $\lambda$ in NSIHT satisfies  the condition  (\ref{rangeoflambda1}),
then the sequence $\{x^n\},$ generated by the NSIHT,  approximates $x_S$ with error
$$ \|x^{p+1} - x_S\|_2 \le \rho^p \|x^0 - x_S\|_2 + \tau \|A x_{\overline{S}} + e\|_2,  $$
where
$$ \rho = \sqrt{3}\left (\delta_{3k} + \sigma_1^2-\frac{\lambda\sigma_1^2}{\epsilon+\sigma_1^2}\right) <1$$
and
$$    \tau = \frac{\sqrt{3} \lambda \sigma_1}{(\epsilon + \sigma_1^2)(1-\rho)}. $$
\end{Thm}

\emph{Proof.}
By the structure of the NSIHT algorithm, we have
\[ u^p = x^p + \lambda (A^TA+\epsilon I)^{-1} A^T(y-Ax^p). \]
In the noisy situation, the measurements of $x$  are given as $y=Ax+e$, where $e \in \mathbb{R}^m$ is the noise vector. By the structure of the algorithm NSIHT, $x^{p+1}= {\cal H}_k (u^p).$ Note that $x_S$ is a $k$-sparse vector where $S=\textrm{supp} ({\cal H}_k(x)). $
By Lemma \ref{lem33}, we get
\begin{equation} \label{dierci}
    \|x^{p+1} - x_S\|_2
    \le \sqrt{3} \| ( u^p - x_S )_{ S^{p+1} \cup S }\|_2,
\end{equation}
where  $ S^{p+1} = \textrm{supp} (x^{p+1}).$ Note that $ y=  A x_S +e'$ where $ e'= Ax_{\overline{S}} +e. $
Thus,
  \begin{align*}
u^p - x_S = & x^p - x_S + \lambda (A^TA + \epsilon I)^{-1} A^T (Ax_S+e' -Ax^p) \\
= & \left( I - \lambda (A^TA + \epsilon I)^{-1} A^T A \right) (x^p - x_S) \\
 &   + \lambda (A^TA + \epsilon I)^{-1} A^T e' .
\end{align*}
Furthermore, we have
 \begin{align} \label{shizi6}
   &   \| (u^p - x_S)_{ S^{p+1} \cup S } \|_2 \nonumber  \\  & =   \|\left[ \left( I - \lambda (A^TA + \epsilon I)^{-1} A^T A \right) (x^p - x_S) \right]_{ S^{p+1} \cup S }  \nonumber \\
     & ~~~~ +\lambda \left(  (A^TA + \epsilon I)^{-1} A^T e' \right)_{ S^{p+1} \cup S } \|_2  \nonumber \\
    & \le   \|\left[ \left( I - \lambda (A^TA + \epsilon I)^{-1} A^T A \right) (x^p - x_S) \right]_{ S^{p+1} \cup S }\|_2 \nonumber  \\
    & ~~~~ + \lambda \| (A^TA + \epsilon I)^{-1} A^T e'  \|_2  \nonumber \\
  &  \leq      (\delta_{3k} + \sigma_1^2-\frac{\lambda\sigma_1^2}{\epsilon+\sigma_1^2})
\|x^p - x\|_2  + \frac{\sigma_1}{\epsilon + \sigma_1^2} \|e'\|_2 ,
\end{align}
where the last inequality follows from Lemma \ref{lemma2} and the bound (\ref{xingxing2}).
Substituting  (\ref{shizi6}) into  (\ref{dierci}) leads to
\begin{align}\label{htp2}
   &  \|x^{p+1} - x_S\|_2 \nonumber  \\
      & \le     \sqrt{3} \Big{ [  } (\delta_{3k} + \sigma_1^2 - \frac{\lambda\sigma_1^2}{\epsilon+\sigma_1^2} )
    \|x^p - x_S\|_2   + \frac{\lambda \sigma_1}{\epsilon + \sigma_1^2} \|e'\|_2  \Big{] } \nonumber \\
& = \rho  \|x^p - x_S\|_2   + \frac{\lambda \sigma_1}{\epsilon + \sigma_1^2} \|e'\|_2  ,
\end{align} where
$$ \rho = \sqrt{3}\left (\delta_{3k} + \sigma_1^2 - \frac{\lambda\sigma_1^2}{\epsilon+\sigma_1^2}\right).
$$
As shown in the proof of Theorem \ref{Thm1},  under the same assumption on $ A$  and  the choice of $ \epsilon$ and $ \lambda,$ we can show that  the  constant $\rho$ above  is smaller than 1. Therefore,
by induction, we obtain the following recovery error bound:
\[ \|x^{p+1} - x_S\|_2 \le \rho^p \|x^0 - x_S\|_2 + \tau \| e' \|_2 , \]
where
 \[
\tau = \frac{\sqrt{3} \lambda \sigma_1}{(\epsilon + \sigma_1^2)(1-\rho)}. \]
 This concludes the proof of the theorem.  ~~  $ \Box. $

The above theorem claims that when the measurements of the signal are slightly inaccurate and  the signal can be sparsely approximated, the recovery of the major information of the signal can be achieved by the algorithm NSIHT.

We now establish the recovery error bound in noisy situations for  the  NSHTP which is a combination of the NSIHT  and the pursuit step. To his purpose, we first state a property of the pursuit step which can be found  in \cite{BFH16} and \cite {zhao2019optimal}.

\begin{Lem} \label{lemmahtp} (\cite {BFH16, zhao2019optimal})
Let $y: =A\widehat{x} +\nu $ be the noisy measurements of the $k$-sparse signal $\widehat{x},$ and let $ v $ be an arbitrarily given  $k$-sparse vector. Then the solution of the pursuit step
$$ z^* = \arg\min_{z} \{\|y-Az\|_2^2:    \textrm{supp} (z) \subseteq  \textrm{supp} (v) \} $$ satisfies that
$$ \|z^*-\widehat{x}\|_2 \leq \frac{1}{\sqrt{1-(\delta_{2k})^2}} \|\widehat{x} -v\|_2 + \frac{\sqrt{1 + \delta_k} } {1 - \delta_{2k}} \| \nu \|_2 . $$
\end{Lem}
Using this property, we now prove the main result for NSHTP.

\begin{Thm} \label{thm3}
Let $\sigma_1$ and $\sigma_m$  be the largest and smallest singular values of  the matrix $A\in \mathbb{R}^{m \times n}$ with $m\ll n.$
Suppose that the  restricted isometry constant of  $A$ satisfies that
$ \delta_{3k} < \frac{1}{2} . $
 Let $y \coloneqq Ax+e$ be the measurements of the signal $x\in\mathbb{R}^n. $
Let  $\epsilon$ be a given parameter satisfying
 \begin{equation} \label{cond-eps} \epsilon > \max\left\{ \sigma_1^2,  ~ \left( \frac{\sigma_1^2-\sigma_m^2}{\frac{1}{2}-\delta_{3k}} - 1 \right) \sigma_1^2\right\}.
\end{equation}
If the stepsize  $\lambda$ in NSHTP is chosen such that
$$  \epsilon + \sigma_1^2 - (\frac{1}{2}-\delta_{3k}  ) \frac{\epsilon + \sigma_1^2}{\sigma_1^2}  < \lambda \leq  \epsilon+\sigma_m^2  ,$$
then the sequence $\{x^n\},$ generated by the NSHTP, satisfies
\begin{equation} \label{errorbound} \| x^{p+1} - x_S\|_2  \leq \rho^p \norm{x^0 - x_S}_2 + \tau \norm{A x_{\overline{S}}+e}_2,
 \end{equation}
where
\[ \rho = \frac{\sqrt{3}\left(\delta_{3k} + \sigma_1^2 - \frac{\lambda\sigma_1^2}{\epsilon+\sigma_1^2}\right)}{\sqrt{1-(\delta_{2k})^2}} < 1,  \]
and
\[
\tau = \frac{1}{1-\rho} \left( \frac{\sqrt{3}\lambda \sigma_1}{\sqrt{1-(\delta_{2k})^2}(\epsilon + \sigma_1^2)} + \frac{\sqrt{1 + \delta_k} } {1 - \delta_{2k}} \right) . \]
\end{Thm}

\emph{Proof.}
We still write $ y= Ax+e$ as  $y = Ax_S+e'$ with $e' = A x_{\overline{S}} + e$, where $S = \text{supp}({\cal H}_k (x))$ and $\overline{S}$ is the complement set of $S$. Note that the intermediate point $\overline{x}^p $ in NSHTP was generated by the NSIHT.
 Based on a result for NSIHT, i.e., the bound (\ref{htp2}), we immediately obtain the following inequality:
\begin{align} \label{nshtpone}
    \| \overline{x}^p- x_S   \|_2 \le  & \sqrt{3}  \Big{ [ } (\delta_{3k} + \sigma_1^2 - \frac{\lambda\sigma_1^2}{\epsilon+\sigma_1^2} )
\|x^p - x_S\|_2 \nonumber  \\
 &  + \frac{\lambda \sigma_1}{\epsilon + \sigma_1^2} \|e'\|_2 \Big { ]  }.
\end{align}
The coefficient in the first term of the right-hand side of (\ref{nshtpone})  is positive. This is guaranteed by $\lambda \le \epsilon + \sigma_m^2$.
Applying Lemma \ref{lemmahtp} yields
\begin{equation} \label{nshtptwo}
    \| x^{p+1} - x_S\|_2
\leq \frac{1}{\sqrt{1-(\delta_{2k})^2}} \| x_S - \overline{x}^p \|_2 + \frac{\sqrt{1 + \delta_k} } {1 - \delta_{2k}} \|e'\|_2 .
\end{equation}
Combining (\ref{nshtpone}) and (\ref{nshtptwo}) leads to
\begin{equation} \begin{aligned} \label{NSHTPsequence}
     \| x^{p+1} &  -   x_S\|_2 \nonumber  \\
& \leq  \frac{\sqrt{3}}{\sqrt{1-(\delta_{2k})^2}} \Big{[} (\delta_{3k} + \sigma_1^2 - \frac{\lambda\sigma_1^2}{\epsilon+\sigma_1^2})
\|x^p - x_S\|_2  \nonumber \\
 &  ~~~~    + \frac{\lambda \sigma_1}{\epsilon + \sigma_1^2} \|e'\|_2  \Big{]} + \frac{\sqrt{1 + \delta_k} } {1 - \delta_{2k}} \|e'\|_2 \nonumber \\
& = \frac{\sqrt{3}(\delta_{3k} + \sigma_1^2 - \frac{\lambda\sigma_1^2}{\epsilon+\sigma_1^2})}{\sqrt{1-(\delta_{2k})^2}} \|x^p - x_S\|_2 \nonumber  \\
 & ~~~~ + \Big{ ( } \frac{\sqrt{3}\lambda \sigma_1}{\sqrt{1-(\delta_{2k})^2}(\epsilon + \sigma_1^2)}
   + \frac{\sqrt{1 + \delta_k} } {1 - \delta_{2k}}   \Big{) }  \|e'\|_2 .
\end{aligned}
\end{equation}
Under the condition of the theorem, we now show that
\begin{equation} \begin{aligned} \label{HTPratio}
    \rho \coloneqq \frac{\sqrt{3}}{\sqrt{1-(\delta_{2k})^2}} \left( \delta_{3k} + \sigma_1^2 - \frac{\lambda\sigma_1^2}{\epsilon+\sigma_1^2} \right) < 1.
\end{aligned} \end{equation}
By the definition of the restricted isometry constant, we see that $ \delta_{2k} \leq \delta_{3k}. $ This together with the condition $\delta_{3k} < 1/2$ implies that
\[ \frac{\sqrt{3}}{2} < \sqrt{1 - (\delta_{3k})^2} \le \sqrt{1 - (\delta_{2k})^2} . \]
Therefore, to ensure that (\ref{HTPratio}) is satisfied, it is sufficient to require that
$
  \delta_{3k} + \sigma_1^2 - \frac{\lambda\sigma_1^2}{\epsilon+\sigma_1^2}  < \frac{1}{2},
$
i.e.,
\[ \lambda > \epsilon + \sigma_1^2 - ( \frac{1}{2}-\delta_{3k} ) \frac{\epsilon + \sigma_1^2}{\sigma_1^2} .  \]
Note that   $\lambda \le \epsilon + \sigma_m^2$ is also required. Therefore the relation (\ref{HTPratio}) is guaranteed provided that the stepsize  $\lambda$ satisfies the following condition:
$$
    \epsilon + \sigma_1^2 - ( \frac{1}{2}-\delta_{3k} ) \frac{\epsilon + \sigma_1^2}{\sigma_1^2} < \lambda \leq  \epsilon+\sigma_m^2  .
$$
In order to guarantee the existence of this range for the stepsize $\lambda, $  it suffices to require that the right-hand side of the above inequality is strictly larger than its left-hand side. This is guaranteed by the choice of
 $\epsilon$ in (\ref{cond-eps}). Therefore the desired  error bound (\ref{errorbound}) is established.  ~~~ $ \Box$ \\

Similar to Theorem \ref{thm2}, the above result claims that when the measurements and the sparse signal are contaminated with noises, the algorithm NSHTP can still recover the signal to a certain level of quality. In the next section, we study the numerical behaviour of the algorithms through random examples of sparse optimization models.

\section{Simulations}

In this section, we demonstrate some simulation results for the NSIHT and NSHTP algorithms.
 All matrices  and sparse vectors are randomly generated. Their entries are assumed to be independent and identically distributed and follow the standard normal distribution.
The supports (i.e., the positions of nonzero entries) of random sparse vectors are chosen according to a uniform distribution.

\subsection{Residual reduction}
 The value of the objective function  in the problem (\ref{PPSS}) is often called the residual.
It was pointed out in \cite{zhao2019optimal}  that the traditional IHT and HTP algorithms may suffer dramatic numerical oscillation in residual reduction during the course of iterations. The simulations on random examples show that the algorithms NSIHT and NSHTP in this paper   can avoid such difficulty in many situations, and thus they can iteratively reduce  the value of the objective of the problem (\ref{PPSS}). Fig. 1 demonstrates such a result and  the residual change during the course of iterations of HTP, NSIHT, and NSHTP. This result was obtained from a  random matrix $A \in \mathbb{R}^{500 \times 1000}$, and a sparse vector $x^* \in \mathbb{R}^{1000}$ with sparsity level $\| x^* \|_0 = 150$ and with exact measurements $y \coloneqq A x^*. $
All algorithms here start from the initial point $x^0=0$. The parameters $\epsilon=1$ and $\lambda=1$ are used in NSIHT and NSHTP for this experiment.
\begin{figure}[H]
\centering
\includegraphics[width=8cm]{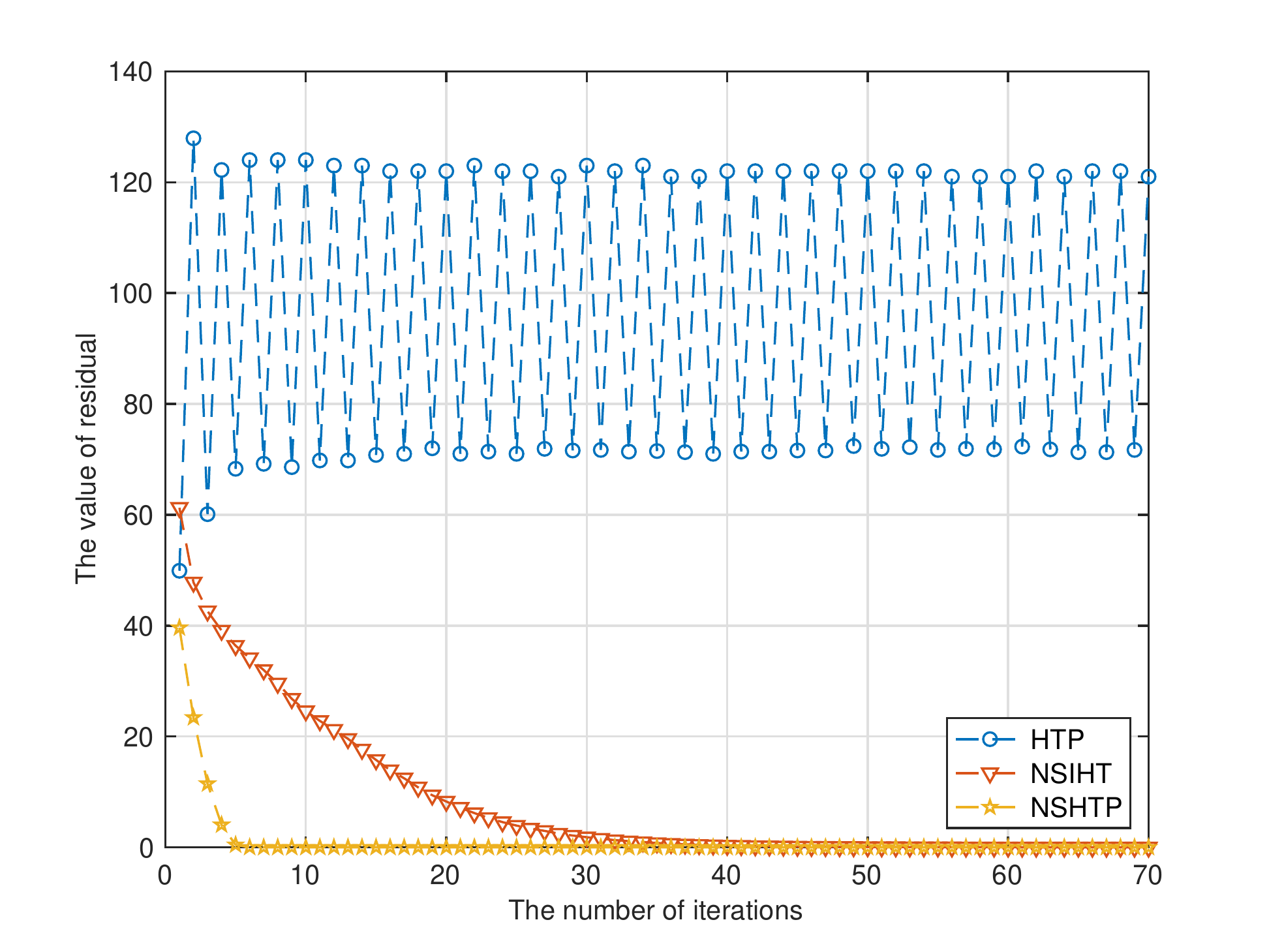}
\caption{The residual reduction in the course of iterations}
\label{graph1}
\end{figure}

Both the NSIHT and NSHTP involve two parameters: $\epsilon$ and $ \lambda.$ Clearly, the choice of these parameters might affect the behaviour of the algorithms.
 By generating random examples   as above (i.e,  $A \in \mathbb{R}^{500 \times 1000}$,  and $x^* \in \mathbb{R}^{1000}$ with   $\| x^* \|_0 = 150$  and $y \coloneqq A x^*  $), we perform NSIHT and NSHTP on such an example to test how the choice of $ (\epsilon, \lambda) $ might  affect the residual reduction in the course of the algorithms.
\begin{figure}
\centering
\subfigure[NSIHT:  $\lambda=1$, different parameter $\epsilon$]{
\includegraphics[width=8cm]{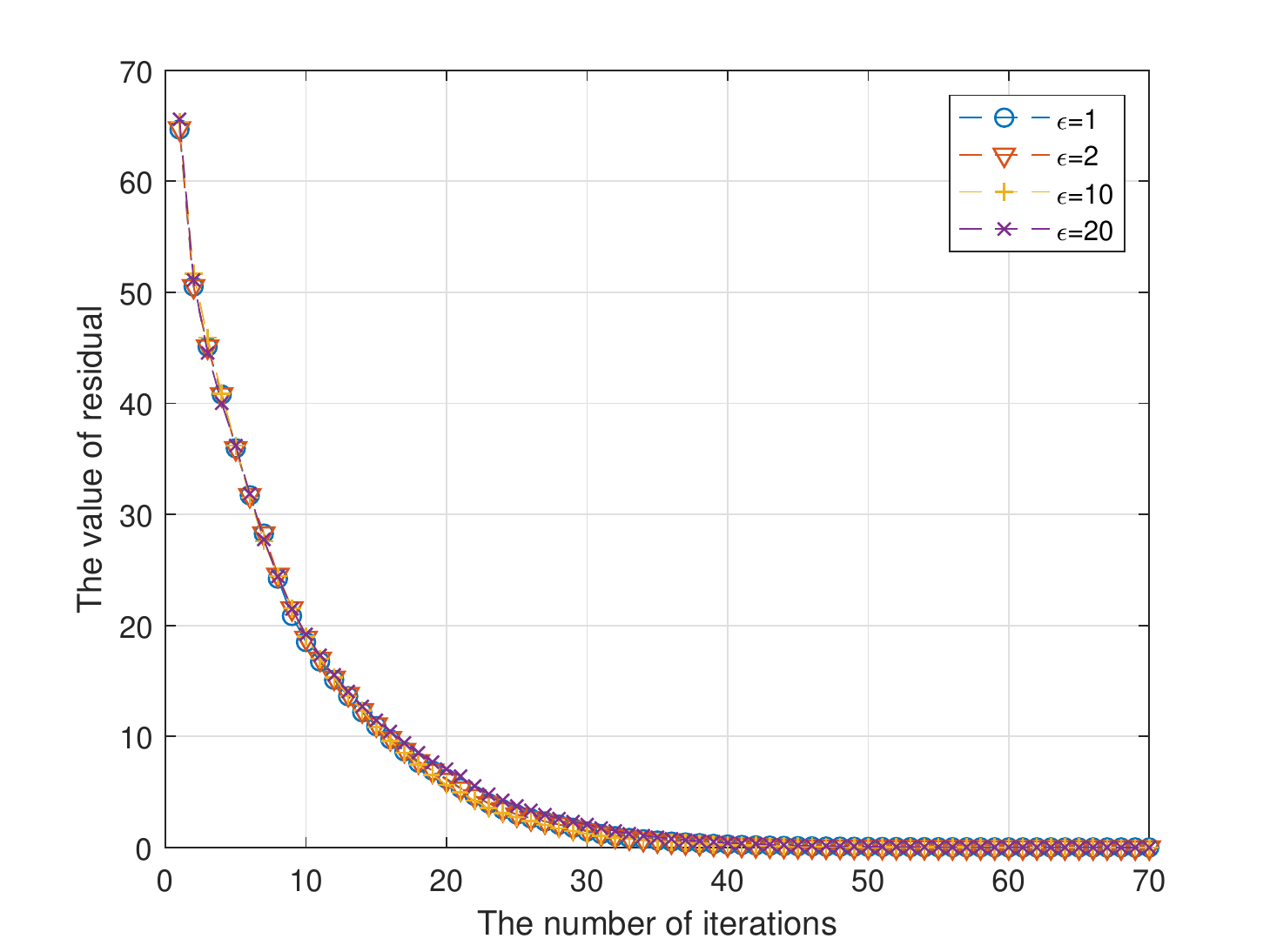} }
\subfigure[NSIHT:  $ \epsilon=1$, different stepsize $\lambda$]{
\includegraphics[width=8cm]{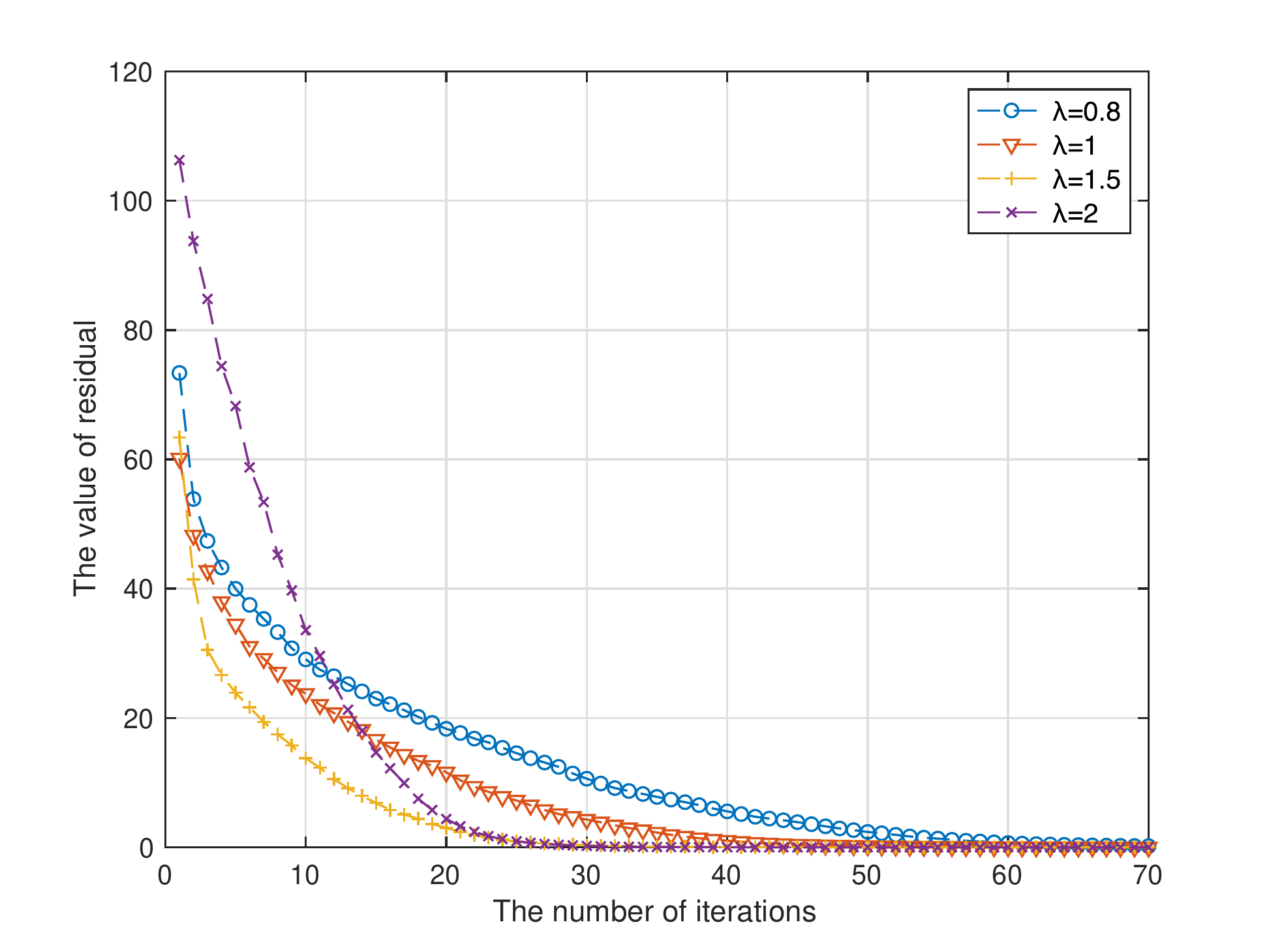} }
\caption{Performance of NSIHT in residual reduction with different $\epsilon$ and $\lambda$}
\label{graph2}
\end{figure}
The result for NSIHT is summarized in Fig. \ref{graph2}, in which (a) is the result for $\lambda=1$ and several different values of $\epsilon,$ and (b) is the result for  $\epsilon=1$ and several different values of $\lambda. $ The simulations indicate that the  performance of the algorithms in residual reduction is closely related to the choice of the stepsize $\lambda$, while  relatively insensitive to the change  of  $\epsilon. $
\begin{figure}
\centering
\subfigure[NSHTP:  $\lambda =1, $ different parameter $\epsilon$]{
\includegraphics[width=8cm]{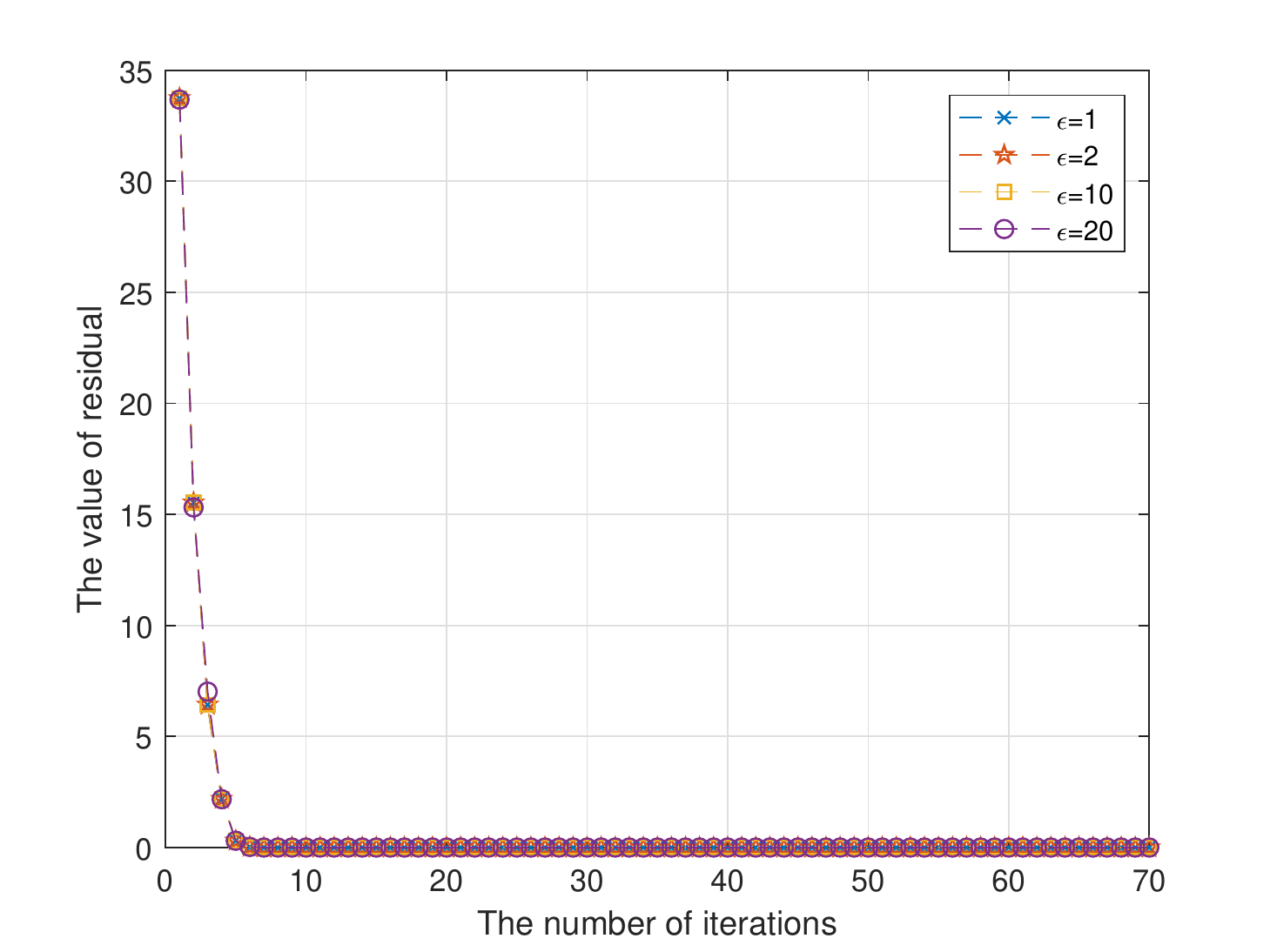} }
\subfigure[NSHTP:  $\epsilon=1, $ different stepsize $\lambda$]{
\includegraphics[width=8cm]{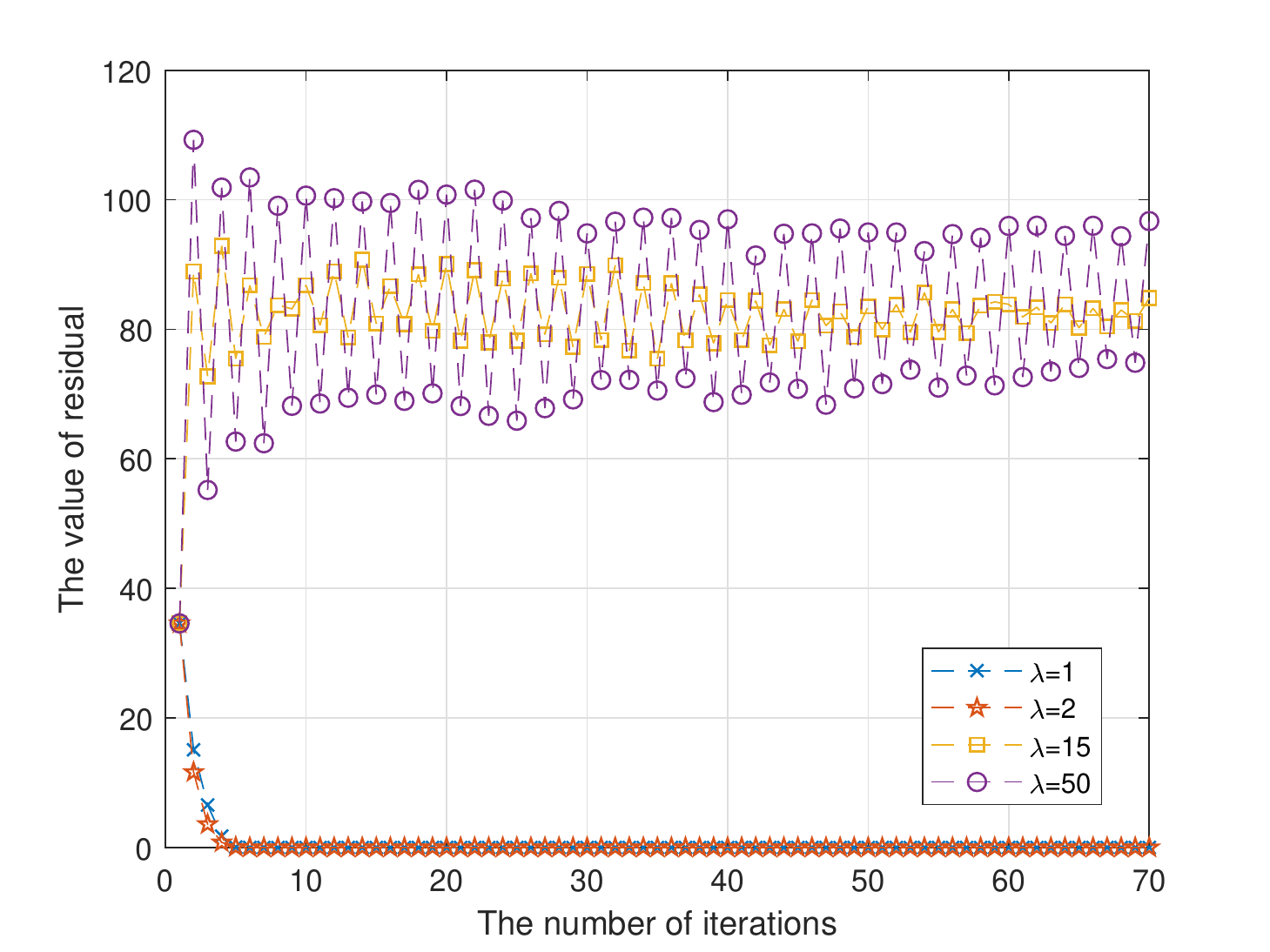} }
\caption{Performance of NSHTP in residual reduction with different $\epsilon$ and $\lambda$}
\label{graph3}
\end{figure}

The results for NSHTP are given in Fig. \ref{graph3}. Similar to NSIHT, the performance  of NSHTP is not very sensitive to the change of  $\epsilon,$  but can be remarkably influenced with the value of the stepsize $ \lambda.$  It can be seen from Fig. \ref{graph3}(b) that the NSHTP  with $ \lambda=1$ or $\lambda=2$  behaves well, however, the oscillation phenomenon in residual reduction were observed when $\lambda=15$ and $\lambda=50 . $

\subsection{Performance of sparse signal recovery}

Several experiments were carried out to demonstrate the performance of the proposed algorithms in sparse signal recovery.
The first one is to test the influence of the ratio between the number of measurements and signal length as well as the number of iterations.
The second experiment is to show how the choice of  the parameters $\lambda$ and  $\epsilon$ might affect the signal recovery performance of the algorithms.
The last one is to examine the recovery capability of the proposed algorithms with noisy or inaccurate measurements of sparse signals.

\subsubsection{Influence of the ratio $m/n$ and number of iterations}

By setting $\epsilon=1$ and $\lambda=1,$  we test the algorithms with three different ratios: $m/n=$ 0.5, 0.4, and 0.3.  The specific random matrices with such ratios in our experiments are taken as  $ 500 \times 1000 $, $ 400 \times 1000 $ and $ 300 \times 1000,  $ respectively.
The algorithms start from the initial point $x^0=0$ and adopt the stopping criterion
\begin{equation} \label{stopcriteria}
    \frac{ \| x^p-x^* \|_2 } {\| x^* \|_2} \le 10 ^{-3}.
\end{equation}

\begin{figure}
\centering
\subfigure[The number of iterations = 20]{
\includegraphics[width=8cm]{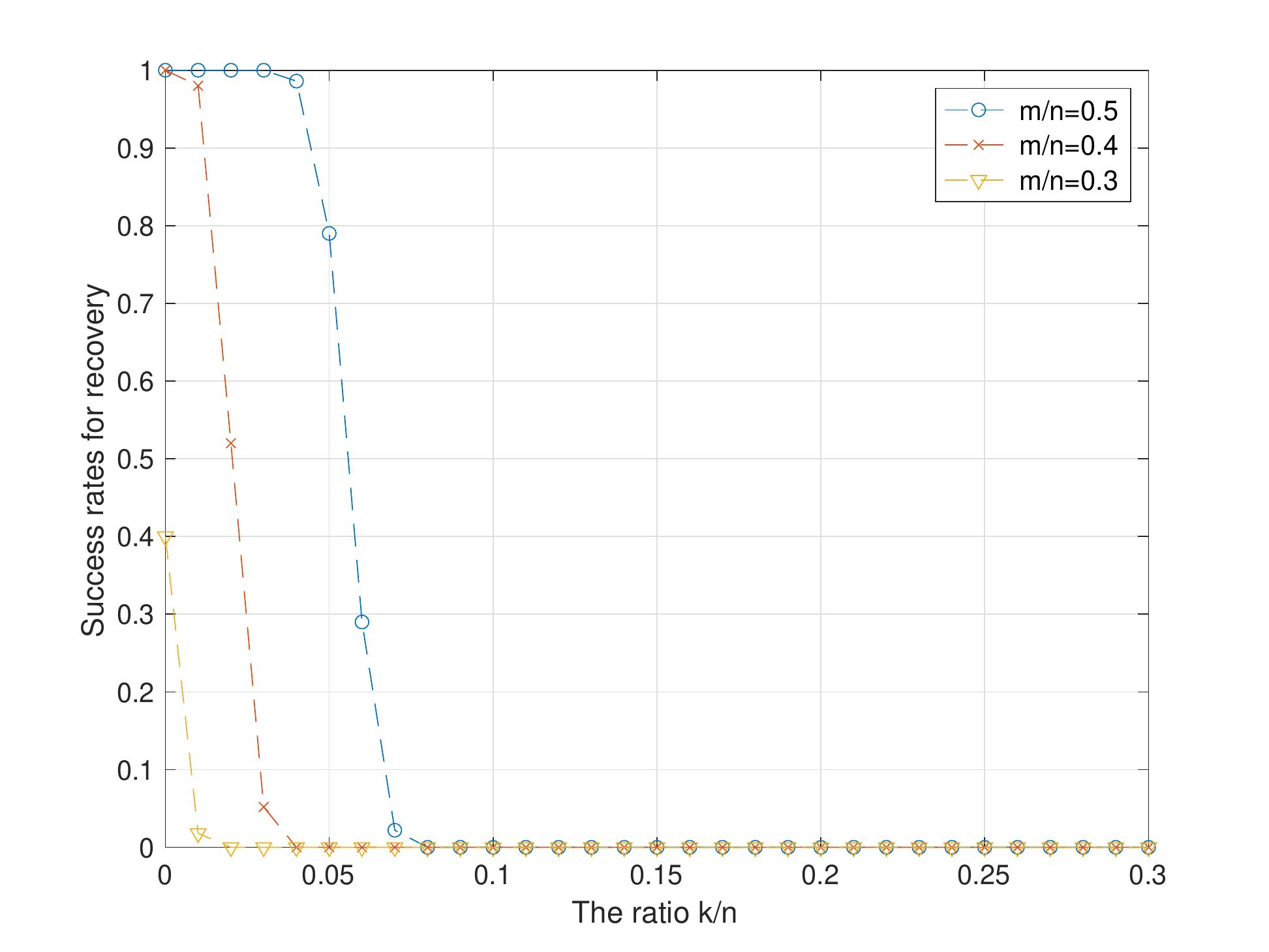}}
\subfigure[The number of iterations = 50]{
\includegraphics[width=8cm]{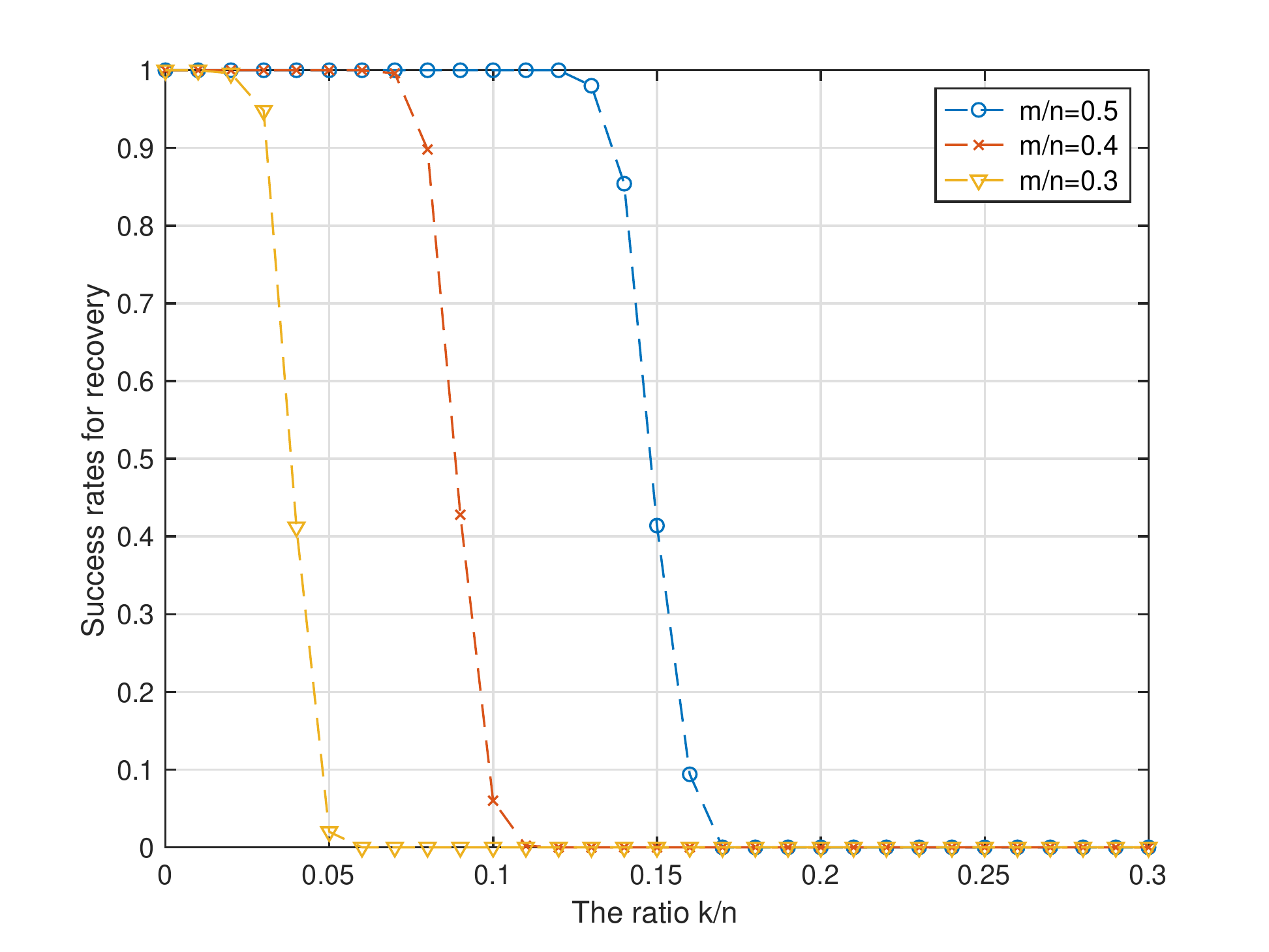}}
\caption{Recovery performance of NSIHT with different number of iterations}
\label{graph4}
\end{figure}

\begin{figure}
\centering
\includegraphics[width=8cm]{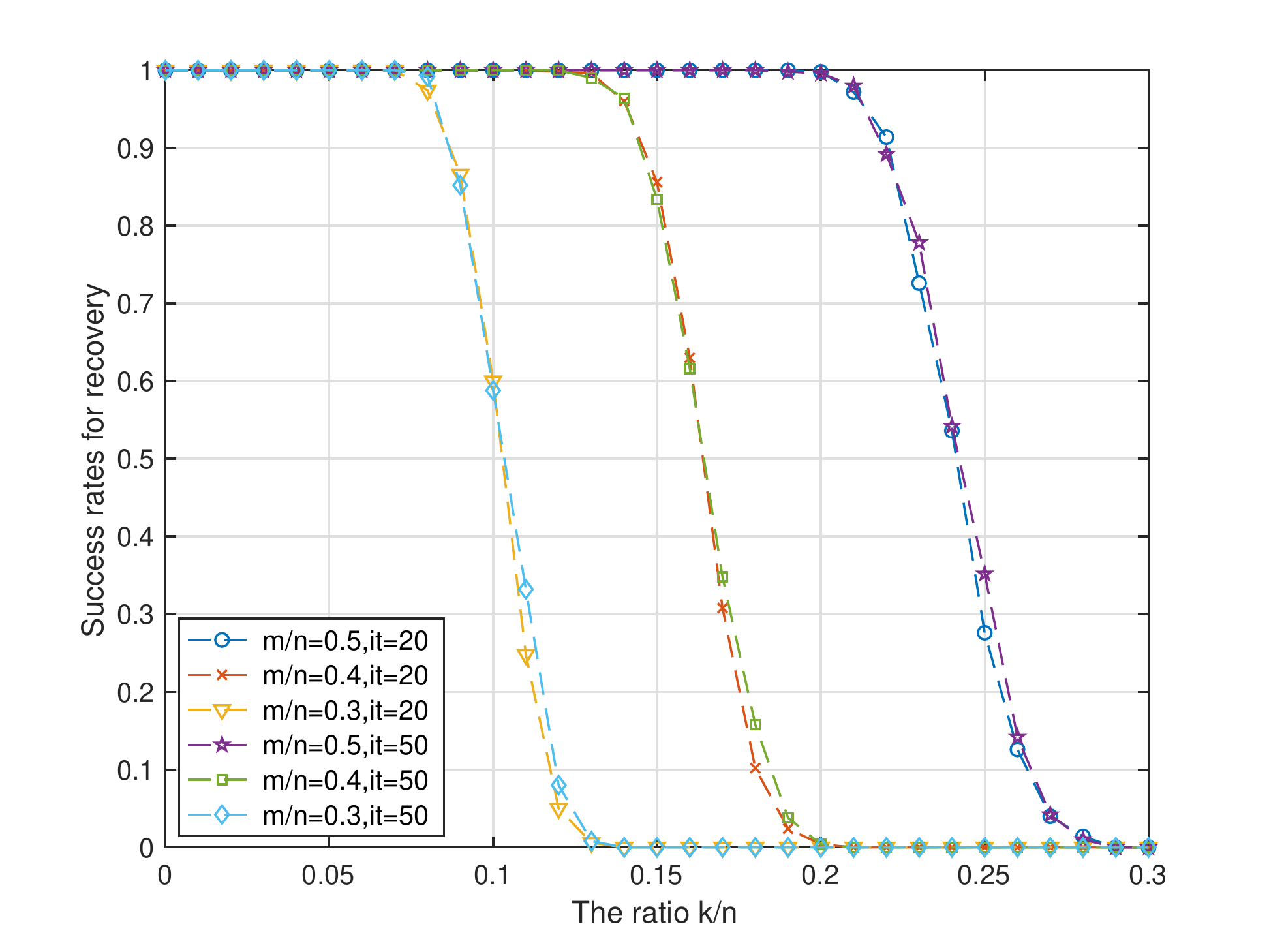}
\caption{Recovery performance of NSHTP   with 20 and 50 iterations}
\label{graph5}
\end{figure}

The results for NSIHT which is performed a total number of 20 and 50 iterations, respectively, are shown in Fig. \ref{graph4}, and the results for NSHTP are shown in Fig. \ref{graph5}.
In these figures, the horizontal axes are the ratios  of the sparsity level $k$ and the signal length $n. $
 The vertical axes record  the success  rate of recovery. The success rate corresponding to each ratio $k/n$ was calculated by performing the algorithms on 250 randomly generated examples of   recovery problems.
From  Fig.  \ref{graph4} and Fig. \ref{graph5},  it can be observed that  performing more iterations generally  produce a better recovery result. Also, we see that the lower the ratio $m/n$ (i.e., the less number of measurements are available),  the narrower the range of sparse signals can be recovered via these algorithms. Compared with Fig. \ref{graph4}, the result in Fig. \ref{graph5} indicates  that in many situations,  using a pursuit step is able to improve the performance of the NSIHT in both convergence speed and signal recovery capability.  When convergent, the NSHTP  usually requires less number of iterations than   NSIHT.

\subsubsection{Influence of the parameters $(\epsilon, \lambda)$ }

\begin{figure}
\centering
\subfigure[NSIHT: 20 iterations]{
\includegraphics[width=8cm]{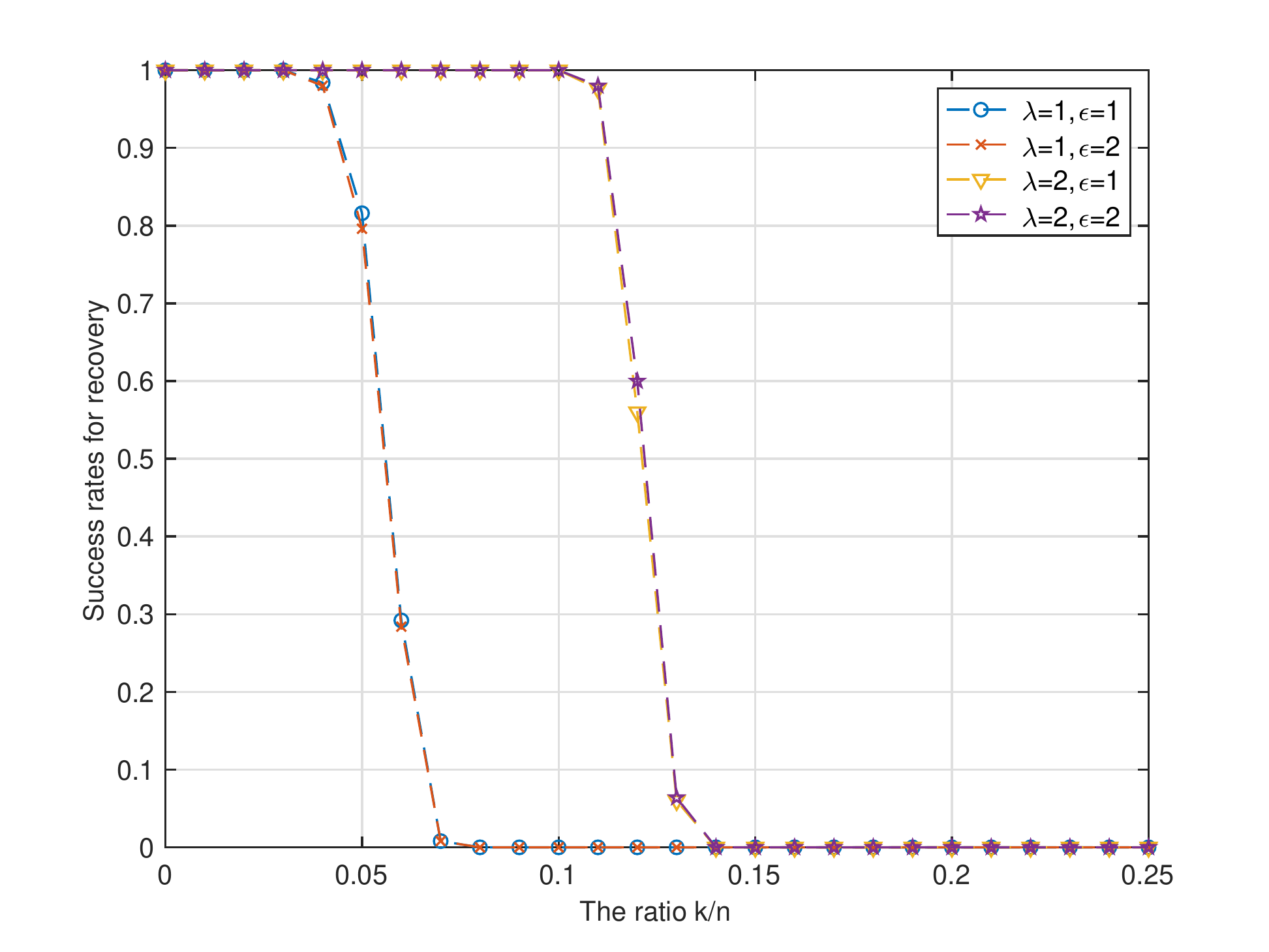}}
\subfigure[NSHTP: 20 iterations]{
\includegraphics[width=8cm]{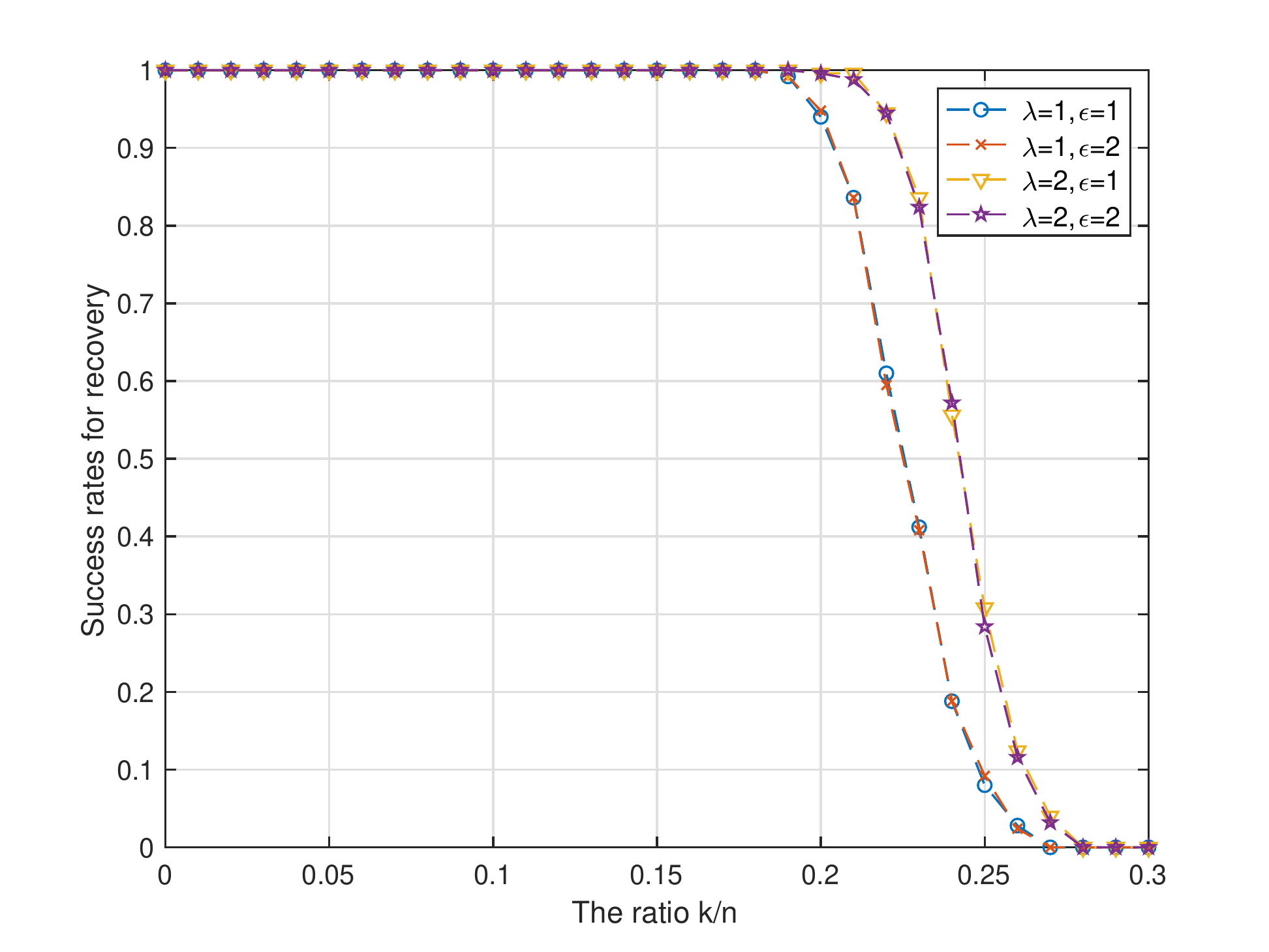}}
\caption{Comparison of NSIHT and NSHTP in signal recovery with different $\lambda$ and $\epsilon$}
\label{graph6}
\end{figure}

\begin{figure}
\centering
\subfigure[NSIHT:  50 iterations]{
\includegraphics[width=8cm]{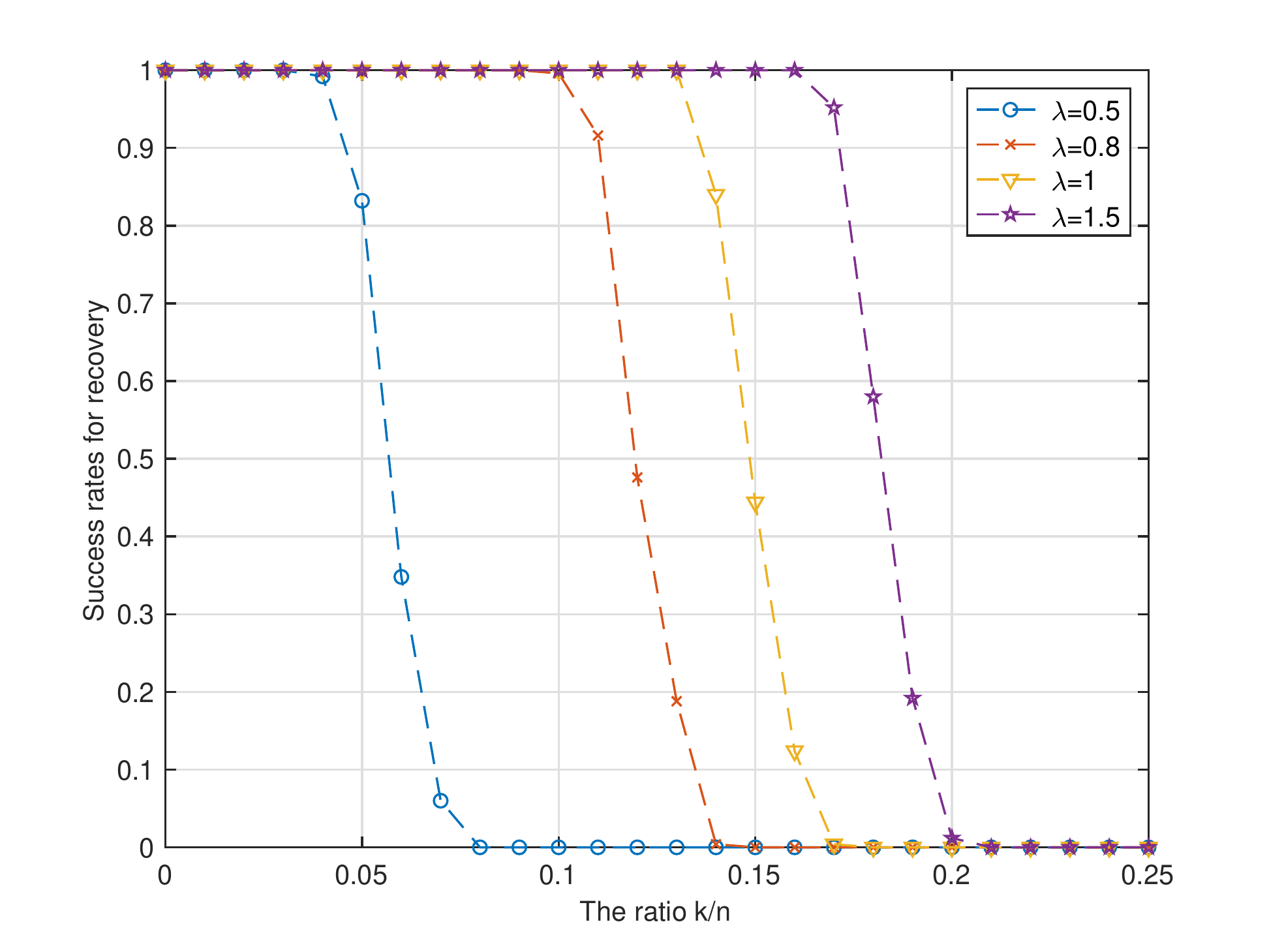}}
\subfigure[NSHTP: 10 iterations]{
\includegraphics[width=8cm]{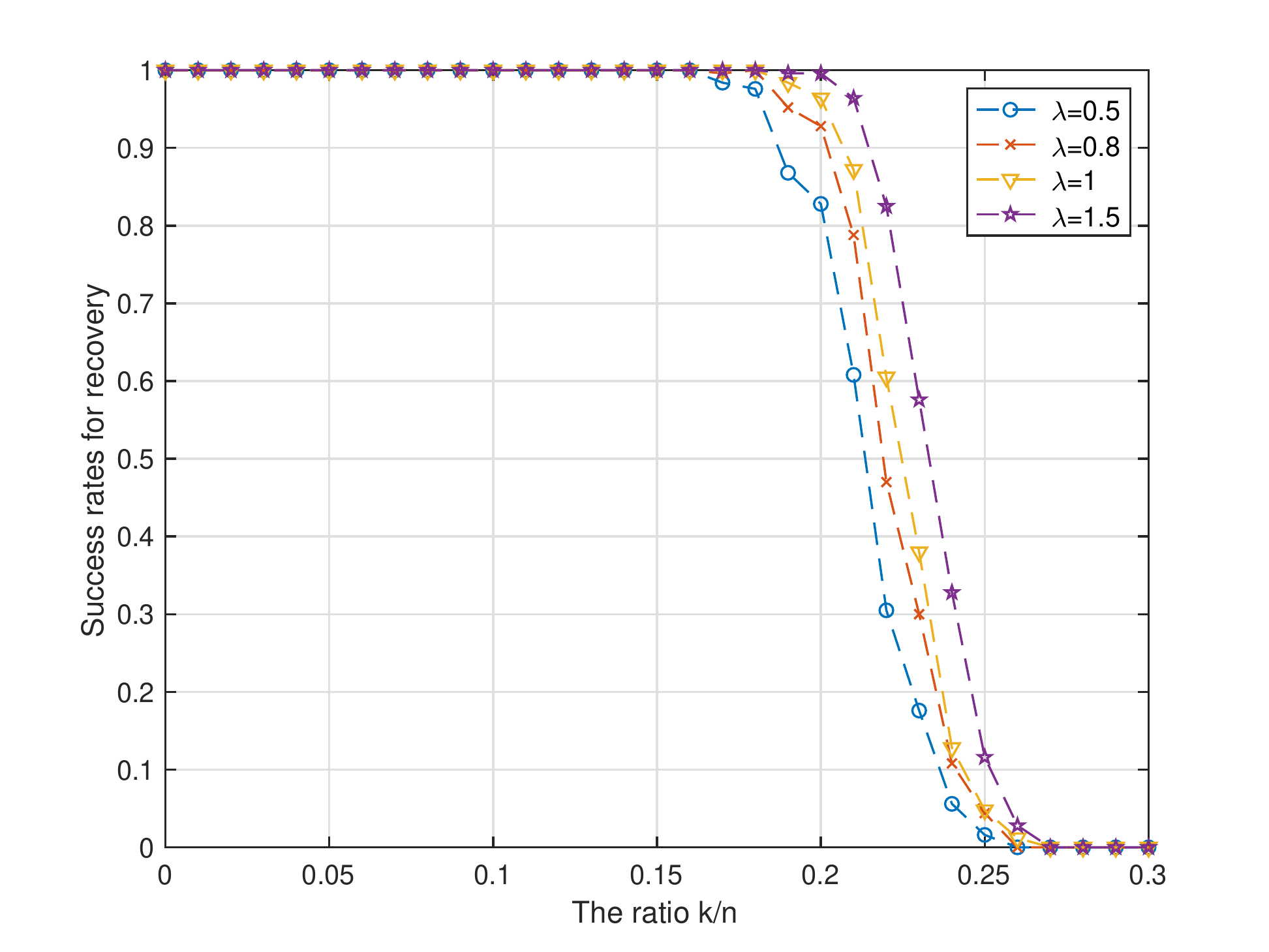}}
\caption{Influence of stepsizes in NSIHT and NSHTP }
\label{graph7}
\end{figure}

 The recovery performance of the algorithms with different choices of their parameters was also examined through simulations. The results for   $\epsilon=1, 2$ and $\lambda=1, 2$, respectively,  are given in Fig. \ref{graph6}. The recovery ability of the algorithms is not very sensitive to the change  of $\epsilon,$ but might be  affected clearly with the change of $\lambda. $ 
Fig. \ref{graph7}  further compares the influence of the stepsize on the recovery ability of NSIHT and NSHTP.
It can be seen that the NSIHT is more sensitive to the change of $\lambda$ than NSHTP.
It seems that $\lambda=1.5, 2 $ are good choices for these algorithms.

\subsubsection{Performance in noisy environments}

\begin{figure}
\centering
\subfigure[The number of iterations = 20]{
\includegraphics[width=8cm]{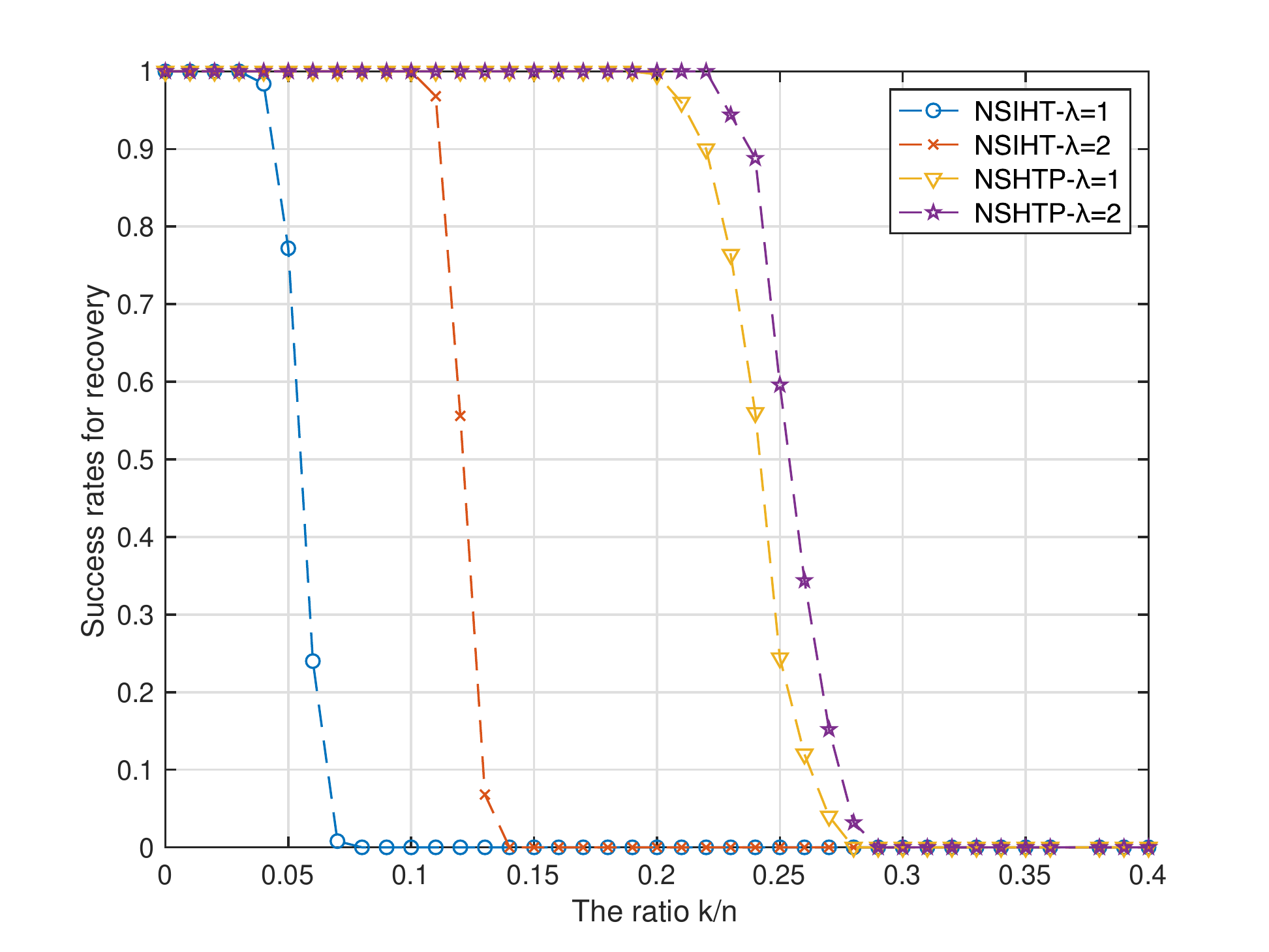}}
\subfigure[The number of iterations = 50]{
\includegraphics[width=8cm]{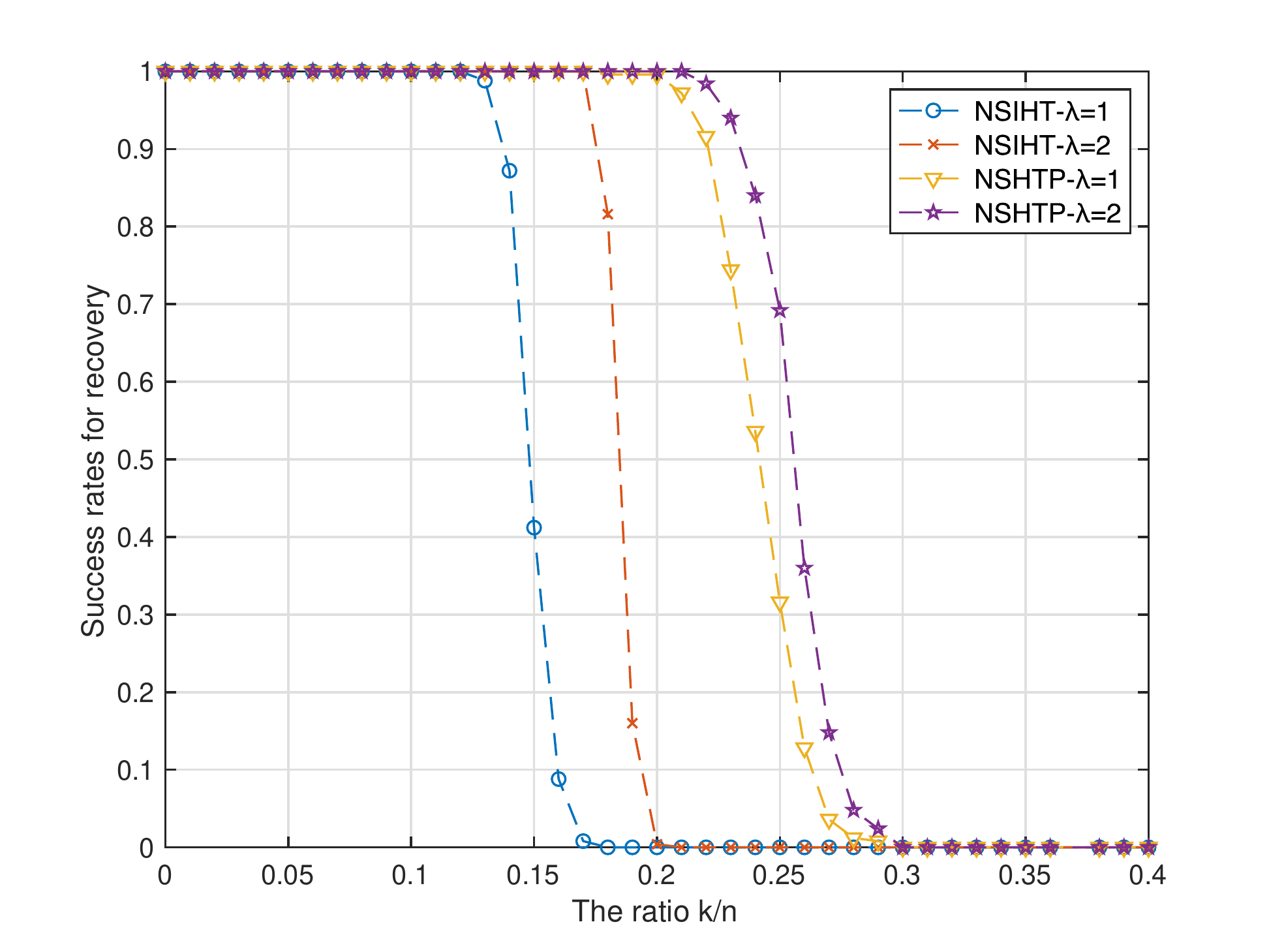}}
\caption{Recovery ability of NSIHT and  NSHTP with inexact measurements}
\label{graph-noisecase}
\end{figure}

Finally, we test the algorithms on  random recovery problems with inaccurate measurements (the measurements are affected by certain  random noises).
In this experiments, the size of random  matrices $A$ is still set as $ 500 \times 1000 $  and the length of the sparse  signal $x^* $  is $n= 1000 .$    The measurements $y = A x^* + e$   is not exact, and the noisy vector $e \in \mathbb{R}^{500}$ is a Gaussian random vector with scale $0.001. $
We still use (\ref{stopcriteria}) as the stopping criterion and the algorithms start from   $x^0=0$.  For each sparsity ratio $k/n$, the algorithms are performed on 250 random examples of the problems. With  $\epsilon = 1, $
the results for two different choices of stepsizes and number of iterations are shown in Fig. \ref{graph-noisecase}.
The simulations indicate  that the NSIHT and NSHTP are able to recover or approximate the signals  from slightly inaccurate measurements of the signals.

\section{Conclusions}

  In this paper, a class of  Newton-step-based hard thresholding algorithms was introduced. The unique feature of the proposed algorithms is that the traditional steepest descent search direction is replaced by a Newton-like search direction. We have proved that with proper choices of the parameter and stepsize, the proposed algorithms can guarantee to recover a sparse signal if the measurement matrix satisfies the standard restricted isometry property. The empirical results indicate that the  new algorithms are efficient for signal recovery and  can also stabilize the objective reduction during the course of iterations.

\end{document}